\documentclass{article}

\usepackage[english]{babel}
\usepackage{appendix}

\usepackage[letterpaper,top=2cm,bottom=2cm,left=3cm,right=3cm,marginparwidth=1.75cm]{geometry}

\usepackage{amsmath}
\usepackage{amssymb}
\usepackage{graphicx}
\usepackage[colorlinks=true, allcolors=blue]{hyperref}
\usepackage{authblk}
\numberwithin{equation}{section}
\newtheorem{theorem}{Theorem}[section]
\newtheorem{lemma}[theorem]{Lemma}
\newtheorem{proposition}[theorem]{Proposition}
\newtheorem{corollary}[theorem]{Corollary}
\newtheorem{definition}[theorem]{Definition}
\newtheorem{remark}[theorem]{Remark}

\newcommand{\scal}{\mathrm{scal}}

\title{Manifolds with positive isotropic curvature of dimension at least nine}
\author[1]{Zhengnan Chen}
\affil[1]{School of Mathematical Sciences, Peking University, Beijing 100871, China}
\date{}
\begin{document}
\maketitle
\footnotetext[1]{E-mail address: \href{mailto:znchen@pku.edu.cn}{znchen@pku.edu.cn}}
\begin{abstract}

In \cite{Brendle19}, Simon Brendle showed that any compact manifold of dimension $n\geq 12$ with positive isotropic curvature and contains no nontrivial incompressible $(n-1)-$dimensional space form is diffeomorphic to a connected sum of finitely many spaces, each of which is a quotient of $S^n$ or $S^{n-1}\times \mathbb{R}$ by standard isometries. We show that this result is actually true for $n\geq9$.
\end{abstract}


\section{Introduction}

\subsection{The positive isotropic curvature condition}

 Positive isotropic curvature condition (PIC for short) was introduced by Micallef and Moore in \cite{MicallefMoore88} in their study of minimal two-spheres in Riemannian manifolds. This curvature condition makes sense in dimension $n\geq4$. In recent years, the study of manifolds with various positive isotropic curvature conditions has become an important subject in the theory of Ricci flow. For example, the conditions PIC1 and PIC2 (variants of the PIC condition) play an important role in the proof of the Differentiable Sphere Theorem \cite{BrendleSchoen07}. We first recall some definitions:

An algebraic curvature operator $R$ is said to be WPIC (weakly positive isotropic curvature) if for any orthonormal four-frame $\{e_1,e_2,e_3,e_4\}$, we have
\begin{equation}
    R_{1313}+R_{1414}+R_{2323}+R_{2424}-2R_{1234}\geq 0.
\end{equation}
$R$ is said to be WPIC1 if for any orthonormal four-frame $\{e_1,e_2,e_3,e_4\}$ and any $\lambda\in [0,1]$,
\begin{equation}
    R_{1313}+\lambda^2R_{1414}+R_{2323}+\lambda^2R_{2424}-2\lambda R_{1234}\geq 0.
\end{equation}
Moreover $R$ is said to be WPIC2 if for any orthonormal 4-frame $\{e_1,e_2,e_3,e_4\}$ and any $\lambda,\mu\in [0,1]$.
\begin{equation}
    R_{1313}+\lambda^2R_{1414}+\mu^2R_{2323}+\lambda^2\mu^2 R_{2424}-2\lambda\mu R_{1234}\geq 0.
\end{equation} 

We call $C_{PIC}, C_{PIC1}, C_{PIC2}$ the cones consisting of algebraic curvature operators that satisfy WPIC, WPIC1 or WPIC2 respectively. Moreover, $R$ is said to be PIC, PIC1 or PIC2 if it's contained in the interior of $C_{PIC1},C_{PIC1}$ or $C_{PIC2}$ respectively. A manifold $M$ is a PIC, PIC1 or PIC2 manifold if there exists a metric $g$ on $M$ such that its curvature operator is PIC, PIC1 or PIC2 respectively. 

Sometimes it is more convenient to use another characterization for the PIC condition. If we extend the operator $R$ complex linearly, then

\begin{enumerate}
    \item $R\in C_{PIC}$ if and only if $R(\varphi,\Bar{\varphi})\geq0$ for all complex two-forms of the form $\varphi=(e_1+ie_2)\wedge (e_3+ie_4)$, where $\{e_1,e_2,e_3,e_4\}$ is an orthonormal four-frame;

    \item $R\in C_{PIC1}$ if and only if $R(\varphi,\Bar{\varphi})\geq0$ for all complex two-forms of the form $\varphi=(e_1+ie_2)\wedge (e_3+i\lambda e_4)$, where $\{e_1,e_2,e_3,e_4\}$ is an orthonormal four-frame and $\lambda\in [0,1]$;

    \item $R\in C_{PIC}$ if and only if $R(\varphi,\Bar{\varphi})\geq0$ for all complex two-forms of the form $\varphi=(e_1+i\mu e_2)\wedge (e_3+i\lambda e_4)$, where $\{e_1,e_2,e_3,e_4\}$ is an orthonormal four-frame and $\lambda,\mu\in [0,1]$;
\end{enumerate}

We remark that PIC1 implies that $R$ has positive Ricci curvature, while PIC2 implies that $R$ has positive sectional curvature.

\subsection{Basic facts regarding the Ricci flow}

Then we recall some basic facts regarding the Ricci flow (one can refer to Chapter 3 of \cite{MorganTian07} for a detailed exposition). As shown in \cite{Hamilton82},\cite{Hamilton86}, under some suitable coordinate frame, the curvature operator $R$ evolves along the Ricci flow by

\begin{equation*}
    \frac{\mathrm{d}}{\mathrm{d}t}R=\triangle R+Q(R),
\end{equation*}
where $Q(R)$ is a quadratic expression of the curvature tensor. More precisely, $Q(R)=R^2+R^{\#}$, where under orthonormal frame we have
\begin{equation*}
    (R^2)_{ijkl}=\sum_{p,q=1}^nR_{ijpq}R_{klpq},
\end{equation*}
and
\begin{equation*}
    (R^{\#})_{ijkl}=2\sum_{p,q=1}^n(R_{ipkq}R_{jplq}-R_{iplq}R_{jpkq}).
\end{equation*}

One important methodology in the theory of Ricci flow is the following: using Hamilton's maximum principle for tensors, we can analyze and control the behavior of the curvature tensor under Ricci flow by establishing a suitable pinching estimate for the Hamiltonian ODE $\frac{\mathrm{d}}{\mathrm{d}t}R=Q(R)$ on the space of algebraic curvature tensors. Combining such pinching estimate with Hamilton's Cheeger-Gromov type compactness theorem (see \cite{Hamilton95ACP} for a reference for the compactness theorem), one deduces that the manifold $(M,g)$ converges to some $(M_{\infty},g_{\infty})$ along the Ricci flow after rescaling. For example, in \cite{Hamilton86}, Hamilton established a pinching estimate for four-dimensional positive curvature operators, hence combining with the compactness theorem, he showed that any four-dimensional compact manifold with positive curvature operator must converge after rescaling to a spherical space form along the Ricci flow, hence must be diffeomorphic to a spherical space form.

\subsection{Böhm-Wilking's pinching cone method}

We then recall the pinching cone method introduced by Christoph Böhm and Burkhard Wilking in \cite{BohmWilking08}, which is used to construct suitable pinching sets and hence to establish suitable pinching estimates for higher dimensional manifolds with various curvature conditions.

Firstly, they introduced the map $l_{a,b}$ on space of algebraic curvature tensors: for any algebraic curvature tensor $R$ and any real numbers $a,b\in \mathbb{R}$, the linear map $l_{a,b}$ is defined to be 

\begin{equation*}
    l_{a,b}(R)=R+bRic(R)\wedge\mathrm{id}+\frac{1}{n}(a-b)\scal(R)\mathrm{id}\wedge\mathrm{id}.
\end{equation*}

Here, the symbol "$\wedge$" denotes the Kulkarni-Nomizu product. If $A$ and $B$ are symmetric bilinear forms on $\mathbb{R}^n$, then \begin{equation*}
    (A\wedge B)_{ijkl}=A_{ik}B_{jl}-A_{il}B_{jk}-A_{jk}B_{il}+A_{jl}B_{ik}.
\end{equation*}

Under the transformation $l_{a,b}$, the scalar curvature of $R$ changes by multiplying a factor $1+2(n-1)a$, and the trace-free Ricci tensor $Ric_0(R)$ changes by multiplying a factor $1+(n-2)b$. 

Now, if we evolves an algebraic curvature operator $R$ by the Hamiltonian ODE $\frac{\mathrm{d}}{\mathrm{d}t}R=Q(R)$, then $S=l_{a,b}^{-1}(R)$ evolves by
\begin{equation*}
    \frac{\mathrm{d}}{\mathrm{d}t}S=Q(S)+D_{a,b}(S),
\end{equation*}
where $D_{a,b}(S)=l_{a,b}^{-1}\circ Q\circ l_{a,b}(S)-Q(S)$ is the difference of $Q(S)$ under the conjugate action of $l_{a,b}$. This property is very important in constructing a suitable family of pinching cones. Here, a family of pinching cones is a family of cones $\{C(b)\},0<b<b_{max}$ for some $b_{max}>0$ that varies continuously and that each $C(b)$ for $0<b<b_{max}$ is a closed, convex, $O(n)-$invariant conic subset of $S_B^2(\mathfrak{so}(n))$, and that each $R\in C(b)\backslash\{0\}$ has positive scalar curvature. Moreover, we require that each $C(b)$ is transversally invariant under the Hamiltonian ODE $\frac{\mathrm{d}}{\mathrm{d}t}R=Q(R)$, which means that for any $R\in \partial  C(b)\backslash\{0\}$, the quadratic $Q(R)$ lies in the interior of the tangent cone $T_{R}C(b)$. 

As computed in \cite{BohmWilking08} (one can also refer to Theorem 11.27 of \cite{Chow2007TheRF} for a detailed exposition), the operator $D_{a,b}(S)$ has the following direct expression:

\begin{equation*}
    \begin{split}
        D_{a,b}(S)=&(2b+(n-2)b^2-2a)Ric_0(S)\wedge Ric_0(S)+2aRic(S)\wedge Ric(S)\\
        &+2b Ric_0(S)^2\wedge\mathrm{id}+\frac{nb^2(1-2b)-2(a-b)(1-2b+nb^2)}{n(1+2(n-1)a)}|Ric_0(S)|^2\mathrm{id}\wedge\mathrm{id}.
    \end{split}
\end{equation*}

In practice, one usually take at first place a continuous family $\{\mathcal{E}(b)\},0<b< b_{max}$ of closed, convex, $O(n)-$ invariant cones for some $b_{max}>0$, and then consider the family $\{C(b)\},0<b< b_{max}$ defined by $C(b)=l_{a,b}(\mathcal{E}(b))$ for suitable $a,b$. To verify that $C(b)$ is transversally invariant under the Hamiltonian ODE is equivalent to verify that $Q(R)$ lies in the interior of $T_{R}C(b)$ for any $R\in \partial C(b)\backslash\{0\}$. To this end, we need to verify that for $S=l_{a,b}^{-1}(R)\in \partial \mathcal{E}(b)\backslash\{0\}$, we have $Q(S)+D_{a,b}(S)$ lies in the interior of $ T_{S}\mathcal{E}(b)$.

If we can find suitable $a,b$ and $\mathcal{E}(b)$ such that $Q(S)+D_{a,b}(S)$ lies in the interior of $ T_{S}\mathcal{E}(b)$ for each $0<b< b_{max}$, then $\{C(b)\},0<b< b_{max}$ is a family of pinching cones. If moreover the family of cones $\{C(b)\},0<b<b_{max}$ converges in the pointed Hausdorff topology to $\mathbb{R}_+\mathrm{id}\wedge\mathrm{id}$ as $b$ approachs to $b_{max}$, then for any compact Riemannian manifold whose curvature operator at each point lies in the interior of $C(0)$, we can construct a suitable pinching set and deduce that the manifold converges to a spherical space form after rescaling (See Theorem 5.1 of \cite{BohmWilking08}).

The pinching cone method is powerful since it helps us to analyze the evolution behavior of the Hamiltonian ODE $\frac{\mathrm{d}}{\mathrm{d}t}R=Q(R)$ in higher dimensions. This method also has many applications. For example, using this method, \cite{BohmWilking08} showed that any compact Riemannian manifold of arbitrary dimension with 2-positive curvature operator (i.e. the sum of two smallest eigenvalues of the curvature operator is positive) is diffeomorphic to a spherical space form.

\subsection{Some related results on various PIC conditions}

Now we look back to the various PIC conditions. it was shown in \cite{BrendleSchoen07} and independently in \cite{Nguyen10} that the WPIC condition is preserved under the Hamiltonian ODE. Moreover, as pointed out in \cite{BrendleSchoen07}, the curvature tensor $R$ of a Riemannian manifold $M$ lies in $C_{PIC1}$ if and only if the curvature tensor $\Tilde{R}$ of $M\times\mathbb{R}$ obtained by extending $R$ to be zero in the additional direction lies in $C_{PIC}$. The curvature tensor $R$ of a Riemannian manifold $M$ lies in $C_{PIC2}$ if and only if the curvature tensor $\hat{R}$ of $M\times\mathbb{R}^2$ obtained by extending $R$ to be zero in the additional two directions lies in $C_{PIC}$. As a consequence, the WPIC1 and WPIC2 condition are also preserved under the Hamiltonian ODE.

In \cite{BrendleSchoen07}, Brendle-Schoen showed that if the sectional curvatures of a Riemannian manifold $(M,g)$ are pointwise strictly $\frac{1}{4}-$pinched, then $(M,g)$ is a PIC2 manifold. Using the pinching cone method, they found suitable pinching estimate and showed that a compact PIC2 manifold must be diffeomorphic to a spherical space form, hence implies the Differential Sphere Conjecture. Later in \cite{Brendle08}, Simon Brendle showed that the above assertion also holds for compact PIC1 manifolds. 

For compact PIC1 manifold, the strategy is still to construct a suitable family of pinching cones. During this process, one crucial ingredient is the positiveness of the Ricci curvature. As pointed out in Proposition 3.2 of \cite{BohmWilking08} (One can also refer to Proposition 11 of \cite{Brendle08}), for any closed, convex, $O(n)-$invariant cone $C$ that is invariant under the Hamiltonian ODE $\frac{\mathrm{d}}{\mathrm{d}t}R=Q(R)$, if any $R\in C\backslash\{0\}$ has positive scalar curvature and nonnegative Ricci curvature, and that $C$ contains all nonnegative curvature operators, then $l_{b+\frac{n-2}{2}b^2,b}(C)$ is transversally invariant under the Hamiltonian ODE for small $b>0$, which is very convenient in use.

For the PIC case, this is more complicated, since we don't have that $Ric\geq0$. The lack of such condition causes us a great difficulty in analyzing and estimating $\frac{\mathrm{d}}{\mathrm{d}t}S$ along the evolution equation $\frac{\mathrm{d}}{\mathrm{d}t}S=Q(S)+D_{a,b}(S)$, and is the main obstacle we face when constructing suitable family of pinching cones. 

In fact, it is not true that a compact PIC manifold must converge to a spherical space form after rescaling. One can easily check that $S^{n-1}\times S^1$ with standard product metric satisfies the PIC condition. A result in \cite{MicallefWang93} shows that any connected sum of finitely many compact spaces, each of which is a quotient of $S^n$ or $S^{n-1}\times \mathbb{R}$ by standard isometries, admits a PIC metric. 

It has been conjectured the converse: is it true a compact PIC manifold must be diffeomorphic to a connected sum of finitely many spaces, each of which is a quotient of $S^n$ or $S^{n-1}\times \mathbb{R}$ by standard isometries?

\subsection{Main result}

Here we give a partially positive answer to the conjecture. The below theorem is the main propose of this paper:

\begin{theorem}\label{thm main}
    Let (M,g) be a compact PIC manifold of dimension $n\geq9$. If $M$ contains no nontrivial incompressible $(n-1)-$dimensional space forms, then $M$ is diffeomorphic to a finite connected sum of spaces, each of which is a quotient of $S^n$ or $S^{n-1}\times \mathbb{R}$ by standard isometries.
\end{theorem}

Our strategy is based on the work of \cite{Brendle19}, where Simon Brendle showed Theorem \ref{thm main} in dimension $n\geq 12$. His proof consists of 3 parts. 

\begin{enumerate}
    \item The first part is to adapt the pinching cone method and establish a suitable pinching estimate for compact PIC manifolds, which ensures that the blow-up limits are uniformly PIC in all dimensions. In dimension $n\geq12$, he showed that the blow-up limits are weakly PIC2.
    \item The second part is to analyze ancient solutions that have bounded curvature, are $\kappa-$noncollapsed, are weakly PIC2, and are uniformly PIC, which is an extension of Perelman's work on characterizing three dimensional $\kappa-$solutions.
    \item  The third part is to develop a surgery theory for Ricci flow starting from a compact PIC manifold, which follows from Perelman's Ricci flow with surgery procedure on three dimensional manifolds.
\end{enumerate}

  We refer the reader to \cite{Perelman02} and \cite{Perelman03} for Perelman's work. We remark that, when $n=4$, Hamilton gave a  pinching estimate for four-dimensional PIC manifolds in \cite{Hamilton97FourmanifoldsWP}, hence showed that the Ricci flow on a four-dimensional compact PIC manifold should only form neck-pinch singularities. Later, \cite{ChenZhu06} adapted and extended Perelman's surgery approach to analyze compact PIC manifolds of dimension four, gave a classification result for four-dimensional compact PIC manifolds that contain no essential incompressible space forms.  By analyzing Ricci flow on orbifolds, Chen-Tang-Zhu in \cite{ChenTangZhu12} gave a complete classification of four-dimensional compact PIC manifolds.

Given Brendle's result for $n\geq12$, it's natural to ask whether we can show similar results for $5\leq n\leq 11$. In fact, condition $n\geq12$ is needed in the first part, while giving a suitable pinching estimate. We notice that the results shown in the second part of \cite{Brendle19} are about characterizing ancient solutions that have bounded curvature, are $\kappa-$noncollapsed, are weakly PIC2, and are uniformly PIC. These results of the second part hold for all $n\geq5$, and do not rely on the pinching estimate in the first part. In the third part, condition $n\geq12$ is needed only because the pinching estimate given in the first part requires that $n\geq12$, and other arguments in the surgery process work for all $n\geq5$. 

As a result, if we can obtain a similar pinching estimate for some lower dimension $n$ with $5\leq n\leq 11$ similar to the pinching estimate in \cite{Brendle19}, then we can adapt the ancient solution analysis and the surgery process in \cite{Brendle19} to classify up to diffeomorphism all compact $n-$dimensional PIC manifolds with no nontrivial incompressible $(n-1)$-dimensional space forms. In this paper, we show that in fact for $9\leq n\leq 11$, we are able to improve the pinching estimate given in \cite{Brendle19} and give a suitable pinching estimate required in the surgery process (see Theorem \ref{thm main cone} and Corollary \ref{cor main pinching estimate} given below), hence we are able to deduce Theorem \ref{thm main}.

To better explain our improvements, we first explain the pinching estimate process in \cite{Brendle19}. This consists of 3 steps.

\begin{enumerate}

\item  Let $K\subseteq int(C_{PIC})$ be a compact set in the interior of $C_{PIC}$. The first step is to construct a family of closed, convex, $O(n)-$invariant pinching sets $\{\mathcal{G}_t^{(0)}: t\in [0,T]\}$ with the property that $\{\mathcal{G}_t^{(0)}:t\in[0,T]\}$ is invariant under $\frac{\mathrm{d}}{\mathrm{d}t}R=Q(R)$. Furthermore, $K\subset\mathcal{G}_0^{(0)}$ and 

    \begin{equation*}
    \begin{split}
        \mathcal{G}_t^{(0)}&\subset \{R:R-\theta \scal\ \mathrm{id}\wedge\mathrm{id}\in C_{PIC}\}\\
        &\cap\{R:Ric_{11}+Ric_{22}-\theta \scal+N\geq 0\}.\\
    \end{split}
    \end{equation*}

The above construction works in dimension $n\geq 5$.

\item Construct a continuous family of pinching cones towards $C_{PIC1}$. This is done by joining 2 families of cones together. One family $\{C(b)\}, 0<b\leq b_{max}$ deforms the cone $\{R\in C_{PIC}:Ric_{11}+Ric_{22}\geq 0\}$ inward, and another family $\{\Tilde{C}(b)\}, 0<b\leq \Tilde{b}_{max}$ deforms the cone $C(b_{max})\cap C_{PIC1}$ outward. 

In dimension $n\geq12$, \cite{Brendle19} showed that the two families can be connected together. Let $\{\hat{C}(b)\}, 0<b< b_{max}+\Tilde{b}_{max}$ denote the connected family, then the family $\{\hat{C}(b)\}$ satisfies the following two properties:

 \begin{enumerate}
        \item For any fixed positive reals $\theta$ and $N$, there exists a small $\beta_0>0$ such that \begin{equation*} \begin{split} (\{R: R-\theta \scal\mathrm{id}\wedge\mathrm{id}\in C_{PIC}\}&\cap\{R: Ric_{11}+Ric_{22}-\theta \scal+N\geq0\} )\\
                &\subset\{R:R+N\mathrm{id}\wedge
                \mathrm{id}\in \hat{C}(\beta_0)\}.
            \end{split}
        \end{equation*}

        \item Let $\{\beta_j\}$ be a increasing sequence such that $\lim_{j\rightarrow +\infty}\beta_j=b_{max}+\Tilde{b}_{max}$, there exists a sequence of positive reals $\{\varepsilon_j\}$ with $\lim_{j\rightarrow +\infty}\varepsilon_j=0$, and
        \begin{equation*}
            \hat{C}(\beta_j)\subset \{R:R+\varepsilon_j\scal\mathrm{id}\wedge\mathrm{id}\in C_{PIC1}\}.
        \end{equation*}
    \end{enumerate}

\item Use the family $\{\hat{C}(b)\},0<b<b_{max}+\Tilde{b}_{max}$ to construct a family of pinching sets $\mathcal{G}_t$ that pinches toward PIC1. The construction use a similar argument as in Theorem 4.1 of \cite{BohmWilking08}. Then, combining the ideas developed in \cite{Brendle08}, one could obtain a pinching set that pinches towards PIC2.
    
\end{enumerate}

 Here we point that, the restriction $n\geq12$ in the work of \cite{Brendle19} appears in the Step 2 above. The arguments in Step 1 and Step 3 actually work for all $n\geq 5$. Hence, if we can use some different constructions and improve the estimates given in \cite{Brendle19} to obtain a sutiable family of pinching cones $\{\hat{C}(b)\}$ for dimension lower than 12, then we can deduce similar results for lower dimensions.

Here is our main improvement: we modify the definition and give more precise estimates for the first family $\{C(b)\},0<b\leq b_{max}$ that pinches $\{R\in C_{PIC}: Ric_{11}+Ric_{22}\geq0\}$ inward. The improved estimation allows us to pinch the first family of pinching cones inside further, and thus we are able to show that the two families of pinching cones $\{C(b)\},0<b\leq b_{max}$ and $\{\Tilde{C}(b)\},0<b\leq \Tilde{b}_{max}$ can be connected together to obtain the desired family $\{\hat{C}(b)\},0<b<b_{max}+\Tilde{b}_{max}$ when the dimension $n$ satisfy $9\leq n\leq 11$. 

In fact, the cone $\{C(b)\},0<b\leq b_{max}$ is obtained by defining $C(b)=l_{a,b}(\mathcal{E}(b))$ for $0<b\leq b_{max}$ and suitable $a,b$ and $\mathcal{E}(b)$. Therefore, in order to show that each cone $\{C(b)\}, 0<b\leq b_{max}$ of the family is transversally invariant under the Hamiltonian ODE, we need to analyze the behavior of $Q(S)+D_{a,b}(S)$ for $S\in\partial\mathcal{E}(b)\backslash\{0\}$. It turns out by calculation that, if one wishes that the term $D_{a,b}(S)$ lies in the interior of $C_{PIC}$, then one needs to show that a suitable lower bound for the sum of two smallest eigenvalues of $Ric(S)$ is preserved under the evolution equation $\frac{\mathrm{d}}{\mathrm{d}t}S=Q(S)+D_{a,b}(S)$. It turns out that to show that this lower bound is preserved, a crucial ingredient is to give an upper bound for $\frac{\mathrm{d}}{\mathrm{d}t}(Ric(S)_{22}-Ric(S)_{11})$, and show that this upper bound is preserved.

Our improvement focuses on the estimation of the two terms $\frac{\mathrm{d}}{\mathrm{d}t}(Ric(S)_{11}+Ric(S)_{22})$ and $\frac{\mathrm{d}}{\mathrm{d}t}(Ric(S)_{11}-Ric(S)_{22})$  along the evolution equation $\frac{\mathrm{d}}{\mathrm{d}t}S=Q(S)+D_{a,b}(S)$. More precisely, when estimating a lower bound for the term $\frac{\mathrm{d}}{\mathrm{d}t}(Ric(S)_{11}+Ric(S)_{22})$, we will encounter a term in the form of $c|Ric_0(S)|^2$, where $c$ is a positive number. It turns out that we could give a lower bound for $|Ric_0(S)|^2$ in terms of $(Ric(S)_{22}-Ric(S)_{11})^2$, which is helpful for our estimation. When estimating an upper bound for the term $\frac{\mathrm{d}}{\mathrm{d}t}(Ric(S)_{11}-Ric(S)_{22})$, we consider the expression for the scalar curvature of $\frac{\mathrm{d}}{\mathrm{d}t}S$ , which is $\frac{\mathrm{d}}{\mathrm{d}t}(\scal(S))=\scal(Q(S))+\scal(D_{a,b}(S))$. It turns out that although the direct expression of $\frac{\mathrm{d}}{\mathrm{d}t}S$ and $Ric(\frac{\mathrm{d}}{\mathrm{d}t}(S))$ are both somewhat complicated, the direct expression of scalar curvature of $\frac{\mathrm{d}}{\mathrm{d}t}S$ has a surprising simplicity. As we will calculate in Section 2, the scalar curvature $\scal(S)$ evolves by

\begin{equation*}
    \begin{split}
        \frac{\mathrm{d}}{\mathrm{d}t}\scal(S)=&\frac{2(1+(n-2)b)^2}{1+2(n-1)a}|Ric(S)|^2+\frac{2((1+2(n-1)a)^2-(1+(n-2)b)^2)}{n(1+2(n-1)a)}(\scal(S))^2\\
        :=&P|Ric(S)|^2+Q(\scal(S))^2.
    \end{split}
\end{equation*}

The expression $1+2(n-1)a$ and $1+(n-2)b$ are just the ratio change of the scalar curvature and the trace-free Ricci tensor of a curvature tensor $R$ under the map $l_{a,b}$. Moreover, we have $P+nQ=2+4(n-1)a$. These properties can simplify our calculation process, and help us to obtain a more precise estimation.

For details, we refer to the proof of Proposition \ref{prop 4th condition for first cone} and Proposition \ref{prop 5th condition of first cone}. In summary, we are able to prove the following.

\begin{theorem}\label{thm main cone}
    Let $9\leq n\leq 11$. Then there exist some $B=B(n)>0$, and a family of closed convex cones $\hat{C}(b), 0<b<B$, such that each cone $\hat{C}(b), 0<b< B$ is transversally invariant under the Hamiltonian ODE $\frac{\mathrm{d}}{\mathrm{d}t}R=Q(R)$, and that
    \begin{enumerate}
        \item For any fixed positive reals $\theta$ and $N$, there exists a small $\beta_0>0$ such that 
        \begin{equation*}
            \begin{split}
                (\{R:R-\theta \scal\mathrm{id}\wedge\mathrm{id}\in C_{PIC}\}&\cap\{R: Ric_{11}+Ric_{22}-\theta \scal+N\geq0\} )\\
                &\subset\{R:R+N\mathrm{id}\wedge
                \mathrm{id}\in \hat{C}(\beta_0)\}.
            \end{split}
        \end{equation*}

        \item Let $\{\beta_j\}$ be a increasing sequence such that $\lim_{j\rightarrow +\infty}\beta_j=B$, there exists a sequence of positive reals $\{\varepsilon_j\}$ with $\lim_{j\rightarrow +\infty}\varepsilon_j=0$, and
        \begin{equation*}
            \hat{C}(\beta_j)\subset \{R:R+\varepsilon_j\scal\mathrm{id}\wedge\mathrm{id}\in C_{PIC1}\}.
        \end{equation*}
    \end{enumerate}
\end{theorem}

As a consequence, we are able to conduct Step 2 of the above pinching estimate in dimension $9\leq n\leq 11$. Now, we follow the arguments of Step 1 and Step 3 as in \cite{Brendle19}, and use the family of pinching cones constructed in Theorem \ref{thm main cone}, we are able to deduce suitable pinching estimate for compact $PIC$ manifolds of dimension $9\leq n\leq 11$. The pinching estimate is formulated as below.

\begin{corollary} \label{cor main pinching estimate}
    (Pinching estimate for PIC manifold, cf. Proposition 1.2 of \cite{Brendle19}) Let $9\leq n\leq 11$. Let $K\subseteq int(C_{PIC})$ be a compact set in the interior of $C_{PIC}$, and $T>0$ be given. Then there exist a small positive real number $\theta$, and a large positive real number $N$, an increasing concave function $f>0$ satisfying $\lim_{s\rightarrow +\infty}\frac{f(s)}{s}=0$, and a continuous family of closed convex, $O(n)-$invariant sets $\{\mathcal{F}_t:t\in[0,T]\}$ such that the family $\{\mathcal{F}_t:t\in[0,T]\}$ is invarinat under Hamilonian ODE $\frac{\mathrm{d}}{\mathrm{d}t}R=Q(R)$, $K\subset \mathcal{F}_0$, and

    \begin{equation*}
    \begin{split}
        \mathcal{F}_t&\subseteq \{R:R-\theta scal \mathrm{id}\wedge \mathrm{id}\in C_{PIC}\}\\
        &\cap\{R:Ric_{11}+Ric_{22}-\theta scal+N\geq 0\}\\
        &\cap\{R:R+f(scal)\mathrm{id}\wedge\mathrm{id}\in C_{PIC2}\}
    \end{split}
    \end{equation*}
for all $t\in[0,T]$.
\end{corollary}

Given this pinching estimate, using Hamilton's maximum principle, we deduce the following pinching estimate for Ricci flow starting from a compact PIC manifolds of dimension $n\geq9$ (cf. Corollary 1.3 of \cite{Brendle19}). That is, for any compact PIC manifold $(M,g_0)$ of dimension $n\geq9$, let $g(t)$ denote the solution to the Ricci flow starting from $g_0$. Then there exist a small positive real number $\theta$, a large positive real number $N$, and an increasing concave function satisfying $\lim_{s\rightarrow+\infty}\frac{f(s)}{s}=0$ such that the curvature tensor of $(M,g(t))$ satisfies $R-\theta\scal\mathrm{id}\wedge\mathrm{id}\in C_{PIC}$, $Ric_{11}+Ric_{22}-\theta\scal+N\geq0$, and $R+f(\scal)\mathrm{id}\wedge\mathrm{id}\in C_{PIC2}$ for all $t\geq0$. Now, following the arguments of the second part and the third part in \cite{Brendle19}, i.e. the ancient solution arguments and the surgery arguments, we are able to deduce Theorem \ref{thm main}. Hence, it suffices to prove Theorem \ref{thm main cone}.

It's worth mentioning that, using our pinching estimate (Theorem \ref{thm main cone}), Hong Huang \cite{HH22}\cite{HH23} claimed that any PIC compact manifold of dimension $n\geq9$ must be diffeomorphic to a connected sum of finitely many spaces, each of which is a quotient of $S^n$ or $S^{n-1}\times\mathbb{R}$ by standard isometries, removing the topological assumption that the manifold contains no incompressible $(n-1)-$dimensional space form.

The rest of this paper is dedicated to the proof of Theorem \ref{thm main cone} and is organized as follows. In Section 2, we give some formulas related to the map $l_{a,b}$. These formulas play a fundamental role in the construction of a desired family of pinching cones. In Section 3, we will give the definition of the first family of pinching cones, and show by computation that each cone in the family is transversally invariant under the Hamiltonian ODE. In Section 4, we will show that the first family of pinching cones we constructed can be connected to the second family constructed in Section 4 of \cite{Brendle19} (with modifications on numerical data), and show that the connected family meets the requirements in Theorem \ref{thm main cone}.

\textbf{Acknowledgements.} The author is grateful to his advisor Professor Gang Tian for encouragements on considering this problem and guidance on an earlier version of this paper. He thanks Duc Nguyen for pointing out a mistake made in a previous version. He also thanks Hong Huang and Ziyi Zhao for helpful discussions.  The author is supported by National Key R\&D Program of China 2020YFA0712800.

\section{Some formulas related to Böhm-Wilking's linear map}

In this section, we give some formulas related to the map $l_{a,b}$ needed in the estimation process of Section 3 and Section 4.

As in the Introduction, for any algebraic curvature operator $R$, $l_{a,b}(R)$ is defined as

\begin{equation}
    l_{a,b}(R)=R+bRic(R)\wedge\mathrm{id}+\frac{1}{n}(a-b)\scal(R)\mathrm{id}\wedge\mathrm{id}.
\end{equation}

Let $D_{a,b}(R)=l_{a,b}^{-1}\circ Q\circ l_{a,b} (R)-Q(R)$, as computed in \cite{BohmWilking08}, we have the formula

\begin{equation*}
    \begin{split}
        D_{a,b}(R) =& (2b+(n-2)b^2-2a)Ric_0(R)\wedge Ric_0(R) \\
        &+ 2a Ric(R)\wedge Ric(R) + 2b^2 Ric_0(R)^2\wedge\mathrm{id}\\
        &+ \frac{nb^2(1-2b)-2(a-b)(1-2b+nb^2)}{n(1+2(n-1)a)}|Ric_0(R)|^2\mathrm{id}\wedge\mathrm{id}.
    \end{split}
\end{equation*}

Tracing the above formula, we obtain the Ricci tensor of $D_{a,b}(R)$

\begin{equation*}
    \begin{split}
        Ric(D_{a,b}(R))=&-4bRic(R)^2+\frac{4}{n}(2b+(n-2)a)\scal(R)Ric(R)\\
        &+2\frac{n^2b^2-2(n-1)(a-b)(1-2b)}{n(1+2(n-1)a)}|Ric_0(R)|^2\mathrm{id}\\
        &+\frac{4}{n^2}(a-b)(\scal(R))^2\mathrm{id}.
    \end{split}
\end{equation*}

Now we calculate $\scal(D_{a,b}(R))$, trace the above formula to get:

\begin{equation*}
    \begin{split}
        \scal(D_{a,b}(R))=& -4b|Ric(R)|^2+\frac{4}{n}(b+(n-1)a)(\scal(R))^2\\
        &+2\frac{n^2b^2-2(n-1)(a-b)(1-2b)}{1+2(n-1)a}|Ric_0(R)|^2
    \end{split}.
\end{equation*}

We recall that if $R$ evolves by the Hamiltonian ODE $\frac{\mathrm{d}}{\mathrm{d}t}R=Q(R)$, then $S=l_{a,b}^{-1}(R)$ evolves by 

\begin{equation*}
    \frac{\mathrm{d}}{\mathrm{d}t}S=Q(S)+D_{a,b}(S).
\end{equation*}

This formula is repeatedly employed to show the transversality of pinching cones to the Hamiltonian ODE. If $S$ evolves by $\frac{\mathrm{d}}{\mathrm{d}t}S=Q(S)+D_{a,b}(S)$, then $Ric(S)$ evolves by

\begin{equation*}
    \begin{split}
        \frac{\mathrm{d}}{\mathrm{d}t}Ric(S)=&2(S\ast Ric(S))-4bRic(S)^2+\frac{4}{n}(2b+(n-2)a)\scal(S)Ric(S)\\
        &+2\frac{n^2b^2-2(n-1)(a-b)(1-2b)}{n(1+2(n-1)a)}|Ric_0(S)|^2\mathrm{id}\\
        &+\frac{4}{n^2}(a-b)(\scal(S))^2\mathrm{id},
    \end{split}
\end{equation*}
where $(S\ast H)_{ik}=\sum_{p,q=1}^nS_{ipkq}H_{pq}$. For the evolution of $\scal(S)$, use $|Ric_0(S)|^2=|Ric(S)|^2-\frac{(\scal(S))^2}{n}$, we deduce that $\scal(S)$ evolves by
\begin{equation*}
    \begin{split}
        \frac{\mathrm{d}}{\mathrm{d}t}\scal(S)=& 2|Ric(S)|^2-4b|Ric(S)|^2+\frac{4}{n}(b+(n-1)a)(\scal(S))^2\\
        &+2\frac{n^2b^2-2(n-1)(a-b)(1-2b)}{1+2(n-1)a}\left(|Ric(S)|^2-\frac{(\scal(S))^2}{n}\right)\\
        =&\frac{2(1+(n-2)b)^2}{1+2(n-1)a}|Ric(S)|^2+\frac{2((1+2(n-1)a)^2-(1+(n-2)b)^2)}{n(1+2(n-1)a)}(\scal(S))^2.
    \end{split}
\end{equation*}

\section{Construction of the first family of pinching cones}\

\subsection{Definitions of the first family of pinching cones}

In this section, we assume $9\leq n\leq11$. Let $b_{max}=\frac{1}{2n+2}$ and $\varepsilon=10^{-3}$, for $0<b\leq b_{max}$, we associate the data $a,\gamma,\rho,\omega,A,P,Q$ by

 \begin{equation}\label{mydata}
    \begin{split}
        a&=a(b)=\frac{(2+(n-2)b)^2}{2(2+(n-3)b)}b,\\
        \gamma&=\gamma(b)=\frac{b}{2+(n-3)b},\\
        \rho&=\rho(b)=b-\frac{2(n-1)\gamma(1-2b)}{n^2}-\frac{2(n-1)(1+\gamma)(n^2b^2-2(n-1)(a-b)(1-2b))}{n^2(1+2(n-1)a)}\\
        \omega&=\omega(b)=(1-\varepsilon)\sqrt{\frac{27(2+(n-2)b)}{8}\frac{b(1+(n-2)b)^2}{n^2\rho^3(2+(n-3)b)^2}},\\
        A&=A(b)=\frac{2+8b}{(n-1)(n-4)}+\frac{4}{n}(2b+(n-2)a),\\
        P&=P(b)=\frac{2(1+(n-2)b)^2}{1+2(n-1)a},\\
        Q&=Q(b)=\frac{2((1+2(n-1)a)^2-(1+(n-2)b)^2)}{n(1+2(n-1)a)}.\\
    \end{split}
\end{equation}

In fact, as calculated in Section 2, the data $P,Q$ satisfies the following: If an algebraic curvature operator $S$ evolves by the ODE $\frac{\mathrm{d}}{\mathrm{d}t}S=Q(S)+D_{a,b}(S)$, then
\begin{equation*}
    \frac{\mathrm{d}}{\mathrm{d}t}(\scal(S))=P|Ric(S)|^2+Q(\scal(S))^2.
\end{equation*}

We define the first family of pinching cones as below. One could compare our definition and Definition 3.1 in \cite{Brendle19}, the fourth condition is different. The advantage of our condition is that we could "push" the cone inside further in the sense that we could let $b_{max}$ bigger. However, to show that our definition makes sense, we need to use a more precise estimate for $\frac{\mathrm{d}}{\mathrm{d}t}(Ric(S)_{22}-Ric(S)_{11})$ (see the proof of Proposition \ref{prop 5th condition of first cone})

\begin{definition}\label{def firstcone}
     For $0<b\leq b_{max}$, let $\mathcal{E}(b)$ denote the set of all algebraic curvature tensors $S$ for which there exists a tensor $T\in S^2(\mathfrak{so}(n))$ satisfying:

     \begin{enumerate}
       \item $T\geq0$,
       \item $S-T\in C_{PIC}$,
       \item $Ric(S)_{11}+Ric(S)_{22}+\frac{2\gamma}{n}scal(S)\geq 0$,
       \item For every orthonormal frame $\{e_1,\cdots,e_n\}$, the following inequality holds:

       \begin{equation}
               Ric(S)_{22}-Ric(S)_{11}\leq\omega^{\frac{1}{2}}(\scal(S))^{\frac{1}{2}}(\sum_{p=3}^n(T_{1p1p}+T_{2p2p}))^{\frac{1}{2}}.
       \end{equation}
     \end{enumerate}

We then define $C(b)=l_{a,b}(\mathcal{E}(b))$.

 \end{definition}

 \begin{remark}
    We write $S^2(\mathfrak{so}(n))$ for the set of curvature-type tensors that do not necessarily satisfy the first Bianchi identity. For $T\in S^2(\mathfrak{so}(n))$, we say $T\geq0$ if $T\in S^2(\mathfrak{so}(n))$ is nonnegative definite. When $T$ satisfies the first Bianchi identity, this coincides with the definition of a non-negative curvature operator. 
    
    Here we abuse a bit the notation for convenience: For $T\in S^2(\mathfrak{so}(n))$, we say $T\in C_{PIC}$ if $T(\phi,\Bar{\phi})\geq0$ for every complex two-form of the form $\phi=(e_1+ie_2)\wedge(e_3+ie_4)$, where $\{e_1,e_2,e_3,e_4\}$ is an orthonormal 4-frame. This is equivalent to say that $T\in C_{PIC}$ if and only if for any orthonormal four-frame $\{e_1,e_2,e_3,e_4\}$ we have

    \begin{equation*}
    T_{1313}+T_{1414}+T_{2323}+T_{2424}+2T_{1342}+2T_{1423}\geq0.
\end{equation*}

Notice that if $T$ satisfies the first Bianchi identity, then the above formula coincides with the definition of $WPIC$ for algebraic curvature operators. 
 \end{remark}

It is clear that $C(b)$ is closed, convex, and $O(n)-$invariant for each $0<b\leq b_{max}$. The following is the main result of this section:

\begin{theorem}\label{thm goal for first family}
    For each $0<b\leq b_{max}$, the cone $C(b)$ is transversally invariant under the Hamiltonian ODE $\frac{\mathrm{d}}{\mathrm{d}t}R=Q(R)$.
\end{theorem}

\subsection{Elementary lemmas}

The rest of this section is dedicated to the proof of Theorem \ref{thm goal for first family}. We first establish some elementary lemmas in need. These results are presented from Lemma \ref{lemincreasing} to Lemma \ref{lem computational results for first family}.

\begin{lemma}\label{lemincreasing}
    For any $0<b\leq b_{max}$ and let $a$ be defined as in (\ref{mydata}), then the 2 functions
    \begin{equation*}
        g(b)=\frac{(1+2(n-2)a)^2}{1+2(n-1)a}\cdot\frac{2+(n-3)b}{1+(n-2)b}
    \end{equation*}
    and
    \begin{equation*}
        h(b)=\frac{(1+2(n-1)a)^2}{(1+2(n-2)a)(1+(n-2)b)^2}
    \end{equation*}
    are both strictly increasing in $b$.
\end{lemma}

\noindent\textit{Proof}. For $g(b)$, direct calculation gives 

\begin{equation*}
    \frac{\mathrm{d}}{\mathrm{d}b}\left(\frac{(1+2(n-2)a)^2}{1+2(n-1)a}\right)=\frac{2(1+2(n-2)a)}{(1+2(n-1)a)^2 }\cdot(2(n-1)(n-2)a+(n-3))\cdot\frac{\mathrm{d}a}{\mathrm{d}b},
\end{equation*}
    
\begin{equation*}
    \frac{\mathrm{d}}{\mathrm{d}b}\left(\frac{2+(n-3)b}{1+(n-2)b}\right)=\frac{1-n}{(1+(n-2)b)^2},
\end{equation*}
use the fact that $\frac{\mathrm{d}a}{\mathrm{d}b}\geq1$, we have

\begin{equation*}
    \begin{split}
        \frac{\mathrm{d}g}{\mathrm{d}b}\geq&\frac{1+2(n-2)a}{(1+2(n-1)a)(1+(n-2)b)}\\
        &\cdot \left(\frac{2(2(n-1)(n-2)a+n-3)(2+(n-3)b)}{1+2(n-1)a}-\frac{(1-n)(1+2(n-2)a)}{1+(n-2)b}\right).
    \end{split}
\end{equation*}

Notice that $a=\frac{(2+(n-2)b)^2}{2(2+(n-3)b)}b<\frac{2+(n+1)b}{2}b\leq\frac{5}{4}b$ and $0<b\leq b_{max}=\frac{1}{2n+2}<\frac{1}{2n}$, we obtain

\begin{equation*}
    \begin{split}
        &\frac{2(2(n-1)(n-2)a+n-3)(2+(n-3)b)}{1+2(n-1)a}-\frac{(n-1)(1+2(n-2)a)}{1+(n-2)b}\\
        >&2(n-3)(2+(n-3)b)-(n-1)(1+2(n-2)a)\\
        >&2(n-3)(2+(n-3)b)-(n-1)(1+\frac{5}{2}(n-2)b)\\
        >&(2(n-3)^2-\frac{5}{2}n^2)b+3n-11\\
        >&-\frac{n(n+12)}{2n}+3n-11\\
        =&\frac{5}{2}n-17>0.
    \end{split}
\end{equation*}

For $h(b)$, it suffices to show that $\frac{1+2(n-1)a}{(1+(n-2)b)^2}$ is strictly increasing in $b$, a direct computation gives

\begin{equation*}
        \frac{\mathrm{d}}{\mathrm{d}b}\left(\frac{1+2(n-1)a}{(1+(n-2)b)^2}\right)=\frac{2}{(1+(n-2)b)^3}\left((n-1)(1+(n-2)b)\frac{\mathrm{d}a}{\mathrm{d}b}-(n-2)(1+2(n-1)a)\right).
\end{equation*}

Notice that $\frac{2+(n-2)b}{2+(n-3)b}$ increases strictly in $b$, hence

\begin{equation*}
    \frac{\mathrm{d}a}{\mathrm{d}b}=\frac{\mathrm{d}}{\mathrm{d}b}\left(\frac{a}{b}\cdot b\right)=\frac{a}{b}+b\frac{\mathrm{d}}{\mathrm{d}b}\left(\frac{a}{b}\right)\geq\frac{a}{b}+b\frac{\mathrm{d}}{\mathrm{d}b}\left(\frac{2+(n-2)b}{2}\right)=\frac{a}{b}+\frac{(n-2)b}{2}.
\end{equation*}

Associating this inequality with the fact that $\frac{a}{b}=\frac{(2+(n-2)b)^2}{2(2+(n-3)b)}>\frac{2+(n-2)b}{2}$ gives us
\begin{equation*}
    \begin{split}
        &(n-1)(1+(n-2)b)\frac{\mathrm{d}a}{\mathrm{d}b}-(n-2)(1+2(n-1)a)\\
        \geq&(n-1)\frac{a}{b}+\frac{(n-1)(n-2)}{2}b+\frac{(n-1)(n-2)^2}{2}b^2-(n-2)-(n-1)(n-2)a\\
        \geq&1+(n-1)(n-2)(b-a)+\frac{(n-1)(n-2)^2}{2}b^2.
    \end{split}
\end{equation*}

Using the fact that $a=\frac{(2+(n-2)^2b)}{2(2+(n-3)b)}b<\frac{2+(n+1)b}{2}\leq\frac{5}{4}b$, one checks directly that $1+(n-1)(n-2)(b-a)+\frac{(n-1)(n-2)^2}{2}b^2>0.$ The assertion follows. $\hfill\square$

\hspace*{\fill}

From the definition of $\rho$ in (\ref{mydata}), one immediately sees that $\rho<b$ and $\lim_{b\rightarrow 0+}\rho=0$. The following lemma gives a characterization of $\rho$ from which we can deduce that $\rho$ is strictly positive for each $0<b\leq b_{max}$

\begin{lemma}\label{lemrho}
    For each $0<b\leq b_{max}=\frac{1}{2n+2}$ and let $\rho$ be defined as in (\ref{mydata}), we have
    \begin{equation}
        \frac{\mathrm{d}\rho}{\mathrm{d}b}>\frac{4}{9}.
    \end{equation}
\end{lemma}

\noindent\textit{Proof}. Notice that 
\begin{equation*}
    \begin{split}
        \frac{\mathrm{d}\rho}{\mathrm{d}b}&=\frac{\mathrm{d}}{\mathrm{d}b}\left(b-\frac{2(n-1)\gamma(1-2b)}{n^2}-\frac{2(n-1)(1+\gamma)(n^2b^2-2(n-1)(a-b)(1-2b))}{n^2(1+2(n-1)a)}\right)\\
        :&=1-P_1-P_2
    \end{split}
\end{equation*}
where
\begin{equation*}
    P_1=\frac{\mathrm{d}}{\mathrm{d}b}\left(\frac{2(n-1)\gamma(1-2b)}{n^2}\right)<\frac{2(n-1)}{n^2}\frac{\mathrm{d}\gamma}{\mathrm{d}b}<\frac{1}{n+1}\leq\frac{1}{10},
\end{equation*}
since the term $1-2b$ is positive and decreasing in $b$ and $\frac{\mathrm{d}\gamma}{\mathrm{d}b}<\frac{1}{2}$. Now, for the term $P_2$, we have

\begin{equation*}
    \begin{split}
        P_2&=\frac{\mathrm{d}}{\mathrm{d}b}\left(\frac{2(n-1)(1+\gamma)(n^2b^2-2(n-1)(a-b)(1-2b))}{n^2(1+2(n-1)a)}\right)\\
        &=\frac{\mathrm{d}}{\mathrm{d}b}\left(\frac{4(n-1)(1+\gamma)((2n-1)b^2+(3n-2)(n-2)b^3+(n-1)(n-2)^2b^4)}{n^2(1+2(n-1)a)(2+(n-3)b)}\right).
    \end{split}
\end{equation*}

Now notice that $1+2(n-1)a$ and $2+(n-3)b$ are both positive and increasing in $b$, combining with $0<b\leq b_{max}=\frac{1}{2n+2}$ and $\frac{\mathrm{d}\gamma}{\mathrm{d}b}<\frac{1}{2}$, we have
\begin{equation*}
    \begin{split}
        P_2\leq&\frac{4(n-1)}{n^2(1+2(n-1)a)(2+(n-3)b)}\frac{\mathrm{d}}{\mathrm{d}b}((1+\gamma)((2n-1)b^2+(3n-2)(n-2)b^3+(n-1)(n-2)^2b^4))\\
        <&\frac{4(n-1)}{n^2(1+2(n-1)a)(2+(n-3)b)}((1+\gamma)(2(2n-1)b+3(3n-2)(n-2)b^2+4(n-1)(n-2)^2b^3))\\
        &+\frac{2(n-1)}{n^2(1+2(n-1)a)(2+(n-3)b)}(2nb^2+3n^2b^3+n^3b^4))\\
        <&\frac{4(n-1)}{n^2(1+2(n-1)b)(2+(n-3)b)}(4nb+9n^2b^2+4n^3b^3)+\frac{n-1}{n^2}(2nb^2+3n^2b^3+n^3b^4)\\
        <&\frac{4(4b+9nb^2+4n^2b^3)}{(1+2(n-1)b))(2+(n-3)b)}+(2b^2+3nb^3+n^2b^4).
    \end{split}
\end{equation*}

We claim that 
\begin{equation*}
    \frac{4(4b+9nb^2+4n^2b^3)}{(1+2(n-1)b))(2+(n-3)b)}<\frac{4}{9}.
\end{equation*}

In fact, the above equation is equivalent to
\begin{equation*}
    36n^2b^3+(81n-2(n-1)(n-3))b^2+(42-5n)b<2,
\end{equation*}
and we verify directly that
\begin{equation*}
    36n^2b^3+(81n-2(n-1)(n-3))b^2+(42-5n)b\leq 36n^2b_{max}^3+(81-2(n-1)(n-3))b_{max}^2<2,
\end{equation*}
hence the claim is true. 

Moreover we can check that
\begin{equation*}
    2b^2+3nb^3+n^2b^4<\frac{1}{100}.
\end{equation*}

Thus 

\begin{equation*}
    \frac{\mathrm{d}\rho}{\mathrm{d}b}=1-P_1-P_2>1-\frac{1}{10}-\frac{4}{9}-\frac{1}{100}>\frac{4}{9},
\end{equation*}
the assertion follows.$\hfill\square$

\begin{lemma}\label{lem computational results for first family}
    Let $9\leq n\leq 11$, let $0<b\leq b_{max}$ and $a,\gamma,\omega,A,P,Q$ be defined as in (\ref{mydata}), the functions $g,h$ be defined as in Lemma \ref{lemincreasing}, then:
    \begin{equation*}
        (1+b\sqrt{n-2})^2< 1+\frac{A}{2}\sqrt{n(n-2)};
    \end{equation*}

\begin{equation*}
        \left(1+\frac{A}{2}\sqrt{n(n-2)}\right)^2\frac{1}{P+nQ}<\frac{\omega}{4},
    \end{equation*}
and
    \begin{equation*}
         \left(1+\frac{A}{2}\sqrt{n(n-2)}\right)\frac{1}{P}<\frac{\omega}{4}.
    \end{equation*}

\end{lemma}

\noindent\textit{Proof.} For the first statement: use the inequality $a>b$, it's clear that

\begin{equation*}
    \begin{split}
        1+\frac{A}{2}\sqrt{n(n-2)}&>1+2b\sqrt{n(n-2)}>1+2(n-2)b\\
        &>1+2\sqrt{n-2}b+(n-2)b^2=(1+b\sqrt{n-2})^2.
    \end{split}
\end{equation*}

For the second statement: we first calculate that
\begin{equation*}
    P+nQ=2+4(n-1)a,
\end{equation*}
hence the statement is equivalent to the following :

\begin{equation}\label{eq in Lem crucial}
    LHS:=\left(1+\frac{A}{2}\sqrt{n(n-2)}\right)^2\cdot\frac{(2+(n-3)b)}{(1+2(n-1)a)(1+(n-2)b)}\leq \frac{1-\varepsilon}{2}\sqrt{\frac{27b(2+(n-2)b)}{8n^2\rho^3}}=:RHS.
\end{equation}

For LHS, observe that
\begin{equation*}
       \left(1+\frac{A}{2}\sqrt{n(n-2)}\right)^2\leq\left(1+\left(\frac{1+4a}{(n-1)(n-4)}+2a\right)(n-1)\right)^2=\left(\frac{n-3}{n-4}\right)^2(1+2(n-2)a)^2.
\end{equation*}

By Lemma \ref{lemincreasing}, we obtain
\begin{equation*}
    LHS\leq\left(\frac{n-3}{n-4}\right)^2g(b)\leq\left(\frac{n-3}{n-4}\right)^2g(b_{max}).
\end{equation*}

Now we show that RHS is strictly decreasing in $b$. Notice that $\rho<b$ and by Lemma \ref{lemrho} , $\frac{\mathrm{d}\rho}{\mathrm{d}b}>\frac{4}{9}$, we have

\begin{equation*}
    \begin{split}
        \frac{\mathrm{d}}{\mathrm{d}b}\left(\frac{b(2+(n-2)b)}{\rho^3}\right)=&\frac{1}{\rho^4}\left((2+2(n-2)b)\rho-3(2+(n-2)b)b\frac{\mathrm{d}\rho}{\mathrm{d}b}\right)\\
        <&\frac{b}{\rho^4}\left((2+2(n-2)b)-\frac{4}{3}(2+(n-2)b)\right)=\frac{b}{\rho^4}\left(-\frac{2}{3}+\frac{2}{3}(n-2)b\right)<0.
    \end{split}
\end{equation*}

We can now verify directly that, for $9\leq n\leq 11$ and $b_{max}=\frac{1}{2n+2}$,  

\begin{equation*}
    LHS\leq \left(\frac{n-3}{n-4}\right)^2g(b_{max})<\frac{1-\varepsilon}{2}\sqrt{\frac{27b_{max}(2+(n-2)b_{max})}{8n^2(\rho(b_{max}))^3}}\leq RHS.
\end{equation*}

The third statement is equivalent to the following:
\begin{equation*}
    LHS:=\left(1+\frac{A}{2}\sqrt{n(n-2)}\right)\cdot\frac{(1+2(n-1)a)(2+(n-3)b)}{(1+(n-2)b)^3}\leq \frac{1-\varepsilon}{2}\sqrt{\frac{27b(2+(n-2)b)}{8n^2\rho^3}}=:RHS.
\end{equation*}

As above, RHS is strictly decreasing in $b$. We now analyze LHS, similarly we have
\begin{equation*}
       1+\frac{A}{2}\sqrt{n(n-2)}\leq\frac{(n-3)(1+2(n-2)a)}{n-4}.
\end{equation*}

Hence by Lemma \ref{lemincreasing},
\begin{equation*}
    LHS\leq \left(\frac{n-3}{n-4}\right)g(b)h(b)\leq \left(\frac{n-3}{n-4}\right)g(b_{max})h(b_{max}).
\end{equation*}

It can be verified directly that, for $9\leq n\leq 11$ and $b_{max}=\frac{1}{2n+2}$,

\begin{equation*}
    LHS\leq\left(\frac{n-3}{n-4}\right)g(b_{max})h(b_{max})<\frac{1-\varepsilon}{2}\sqrt{\frac{27b_{max}(2+(n-2)b_{max})}{8n^2(\rho(b_{max}))^3}}\leq RHS,
\end{equation*}
the assertion follows.$\hfill\square$

\subsection{Proof of Theorem \ref{thm goal for first family}}

After these preparations, we now begin the proof of Theorem \ref{thm goal for first family}. It suffices to show that the cone $\mathcal{E}(b)$ is transversally invariant under the ODE $\frac{\mathrm{d}}{\mathrm{d}t}S=Q(S)+D_{a,b}(S)$ for each $0<b\leq b_{max}$, this is guaranteed by Proposition \ref{prop 2ndprop of first cone} to Proposition \ref{prop 5th condition of first cone}.

The following three propositions are adaptions of work of \cite{Brendle19}. We present the proof here for the convenience of the reader.

\begin{proposition}\label{prop 1stprop of first cone}
 (cf. Lemma 3.6 in \cite{Brendle19})   Suppose $S\in\mathcal{E}(b)\backslash\{0\}$, then $D_{a,b}(S)$ lies in the interior of $C_{PIC}$.
\end{proposition}

\textit{Proof.} Recall that 

\begin{equation*}
    \begin{split}
        D_{a,b}(S) =& (2b+(n-2)b^2-2a)Ric_0(S)\wedge Ric_0(S) \\
        &+ 2a Ric(S)\wedge Ric(S) + 2b^2 Ric_0(S)^2\wedge\mathrm{id}\\
        &+ \frac{nb^2(1-2b)-2(a-b)(1-2b+nb^2)}{n(1+2(n-1)a)}|Ric_0(S)|^2\mathrm{id}\wedge\mathrm{id}.
    \end{split}
\end{equation*}

Firstly, we have

\begin{equation*}
    nb^2(1-2b)-2(a-b)(1-2b+nb^2)=\frac{b^2}{2+(n-3)b}(2+(n-8)b-2(n+2)(n-2)b^2-n(n-2)^2b^3)>0
\end{equation*}
for $0<b\leq b_{max}$. It is left to show that 

\begin{equation*}
    U:=(2b+(n-2)b^2-2a)Ric_0(S)\wedge Ric_0(S)+2aRic(S)\wedge Ric(S)+2b^2Ric_0(S)^2\wedge\mathrm{id}\in C_{PIC}.
\end{equation*}

It suffices to show that for any linearly independent vectors $\zeta,\eta$ satisfying $g(\zeta,\zeta)=g(\zeta,\eta)=g(\eta,\eta)=0$, we have $U(\zeta,\eta,\Bar{\zeta},\Bar{\eta})\geq0$. As $Ric(S)$ is unitarily diagonalizable on $\mathrm{span}_{\mathbb{C}}\{\eta,\zeta\}$, we can find $z,w\in\mathrm{span}_{\mathbb{C}}\{\eta,\zeta\}$ such that $g(z,\Bar{z})=g(w,\Bar{w})=2, g(z,\Bar{w})=0$ and $Ric(S)(z,\Bar{w})=0$. The identities $g(\zeta,\zeta)=g(\zeta,\eta)=g(\eta,\eta)=0$ give $g(z,z)=g(w,w)=g(z,w)=0$. Consequently, we may write $z=e_1+ie_2$ and $w=e_3+ie_4$ for some orthonormal four-frame $\{e_1,e_2,e_3,e_4\}\subset \mathbb{R}^n$. Using the identity $Ric(S)(z,\Bar{w})=0$, we obtain

\begin{equation*}
    \begin{split}
        U(z,w,\Bar{z},\Bar{w})=&2(2b+(n-2)b^2-2a)(Ric_0(S)_{11}+Ric_0(S)_{22})(Ric_0(S)_{33}+Ric_0(S)_{44})\\
        &+4a(Ric(S)_{11}+Ric(S)_{22})(Ric(S)_{33}+Ric(S)_{44})\\
        &+4b^2((Ric_0(S)^2)_{11}+(Ric_0(S)^2)_{22}+(Ric_0(S)^2)_{33}+(Ric_0(S)^2)_{44})\\
        \geq&2(2b+(n-2)b^2-2a)(Ric_0(S)_{11}+Ric_0(S)_{22})(Ric_0(S)_{33}+Ric_0(S)_{44})\\
        &+4a(Ric(S)_{11}+Ric(S)_{22})(Ric(S)_{33}+Ric(S)_{44})\\
        &+2b^2((Ric_0(S)_{11}+Ric_0(S)_{22})^2)+(Ric_0(S)_{33}+Ric_0(S)_{44})^2).\\
    \end{split}
\end{equation*}

Let $x:=\frac{n}{\scal(S)}(Ric_0(S)_{11}+Ric_0(S)_{22})$ and $y:=\frac{n}{\scal(S)}(Ric_0(S)_{33}+Ric_0(S)_{44})$. Recall the third condition in Definition \ref{def firstcone}, we have $x,y\geq-2\gamma-2$. Let $D=\{(x,y)\in \mathbb{R}^2:x,y\geq-2\gamma-2\}$, we show that for any $(x,y)\in D$,

\begin{equation*}
    f(x,y)=(2b+(n-2)b^2-2a)xy+2a(x+2)(y+2)+b^2(x^2+y^2)\geq0.
\end{equation*}

In fact, the Hessian of $f$ is given by

\begin{equation*}
    \mathrm{Hess}f=\begin{pmatrix}
        2b^2 & 2b+(n-2)b^2  \\
        2b+(n-2)b^2 & 2b^2
    \end{pmatrix}.
\end{equation*}

It follows that $\mathrm{Hess}f$ has two eigenvalues of opposite signs. Since $f(x,y)\rightarrow+\infty$ as $(x,y)\in D$ goes to infinity, $f$ attains its minimum on $\partial D$. Now, one checks directly that

\begin{equation*}
    f(-2\gamma-2,y)=f(-\frac{2(2+(n-2)b)}{2+(n-3)b},y)=b^2y^2\geq0.
\end{equation*}

As a result, $U(z,w,\Bar{z},\Bar{w})=\frac{2(\scal(S))^2}{n^2}f(x,y)\geq0$. Since $\mathrm{span}_{\mathbb{C}}\{\zeta,\eta\}=\mathrm{span}_{\mathbb{C}}\{z,w\}$, it follows that $U(\zeta,\eta,\Bar{\zeta},\Bar{\eta})\geq0$, and the assertion follows.$\hfill\square$

\hspace*{\fill}

We evolve the tensor $T$ as in \cite{Brendle19}: use Proposition \ref{prop 1stprop of first cone} to find a small positive number $\varepsilon$ (depending on $b$) such that $D_{a,b}(S)-2\varepsilon(\scal(S))^2\mathrm{id}\wedge\mathrm{id}\in C_{PIC}$ for all $S\in \mathcal{E}(b)$. Then, we evolve $T$ by the ODE
\begin{equation*}
    \frac{\mathrm{d}}{\mathrm{d}t}T=S^2+\varepsilon(\scal(S))^2\mathrm{id}\wedge\mathrm{id}.
\end{equation*}

Directly we see that $S^2$ is a nonnegative symmetric operator in $S^2(\mathfrak{so}(n))$, thus we have

\begin{proposition}\label{prop 2ndprop of first cone}
(cf. Lemma 3.7 in \cite{Brendle19})    Suppose $S\in\mathcal{E}(b)\backslash\{0\}$, $T\geq0$ and $S-T\in C_{PIC}$. Then
    \begin{equation*}
        \frac{\mathrm{d}}{\mathrm{d}t}T>0.
    \end{equation*}
\end{proposition}

We also have the following.

\begin{proposition}\label{prop 3rdprop of first cone}
(cf. Lemma 3.8 in \cite{Brendle19})    Suppose $S\in\mathcal{E}(b)\backslash\{0\}$, $T\geq0$ and $S-T\in C_{PIC}$. Then $\frac{\mathrm{d}}{\mathrm{d}t}(S-T)$ lies in the interior of the tangent cone $T_{S-T}C_{PIC}$.
\end{proposition}

\noindent\textit{Proof.} From the evolution equation of $S$ and $T$, we obtain that

\begin{equation*}
    \frac{\mathrm{d}}{\mathrm{d}t}(S-T)=S^{\#}+D_{a,b}(S)-\varepsilon(\mathrm{scal}(S))^2\mathrm{id}\wedge\mathrm{id}.
\end{equation*}

By the definition of $\varepsilon$, it is left to show that $S^{\#}\in T_{S-T}C_{PIC}$. Equivalently, given any complex 2-form $\varphi\in \mathfrak{so}(n,\mathbb{C})$ with the property that $\varphi=(e_{1}+ie_{2})\wedge(e_3+ie_4)$ for some orthonormal four-frame $\{e_1,e_2,e_3,e_4\}$ satisfying $(S-T)(\varphi,\Bar{\varphi})=0$, we show that $S^{\#}(\varphi,\Bar{\varphi})\geq 0.$ 

The ideas are from the proof of Theorem 1 in \cite{Wilking10}. We know that $S^{\#}(\varphi,\Bar{\varphi})=-\frac{1}{2}\mathrm{tr}(\mathrm{ad}_{\varphi}S\mathrm{ad}_{\Bar{\varphi}}S)$. Notice that for any $\theta \in \mathfrak{so}(n,\mathbb{C})$, we have
\begin{equation*}
    (S-T)(\mathrm{Ad}_{\mathrm{exp}(t\theta)}\varphi,\mathrm{Ad}_{\mathrm{exp}(t\Bar{\theta})}\Bar{\varphi})\geq0.
\end{equation*}

Differentiating twice the above equation at $t=0$, and use the fact that $\mathrm{Ad}_{\mathrm{exp}(t\theta)}=\mathrm{exp}(t\mathrm{ad}_{\theta})$, we obtain the following:

\begin{equation*}
    0\leq 2(S-T)(\mathrm{ad}_{\theta}\varphi,\mathrm{ad}_{\Bar{\theta}}\Bar{\varphi})+(S-T)(\mathrm{ad}_{\theta}\mathrm{ad}_{\theta}\varphi,\Bar{\varphi})+(S-T)(\varphi,\mathrm{ad}_{\Bar{\theta}}\mathrm{ad}_{\Bar{\theta}}\Bar{\varphi}).
\end{equation*}

Replace $\theta$ by $i\theta$, we see that the first summand in the above inequality is left unchanged while the other two summands change their sign. Therefore we add the two inequalities to obtain that

\begin{equation*}
    0\leq(S-T)(\mathrm{ad}_{\theta}\varphi,\mathrm{ad}_{\Bar{\theta}}\Bar{\varphi}).
\end{equation*}

Since $T$ is nonnegative, we have $0\leq S(\mathrm{ad}_{\theta}\varphi,\mathrm{ad}_{\Bar{\theta}}\Bar{\varphi})=S(\mathrm{ad}_{\varphi}\theta,\mathrm{ad}_{\Bar{\varphi}}\Bar{\theta})$ for any $\theta\in\mathfrak{so}(n,\mathbb{C})$. This is equivalent to saying that the operator $-\mathrm{ad_{\Bar{\varphi}}}S\mathrm{ad}_{\varphi}$ is non-negative Hermitian on $\mathfrak{so}(n,\mathbb{C})$. Recall that we want to show that $-\mathrm{tr}(\mathrm{ad}_{\varphi}S\mathrm{ad}_{\Bar{\varphi}}S)=-\mathrm{tr}(\mathrm{ad}_{\Bar{\varphi}}S\mathrm{ad}_{\varphi})\geq0$, using the spectral decomposition of $-\mathrm{ad_{\Bar{\varphi}}}S\mathrm{ad}_{\varphi}$, it suffices to show that $S$ is nonnegative on the image of the operator $\mathrm{ad}_{\Bar{\varphi}}S\mathrm{ad}_{\varphi}$. Clearly we have that the image is contained in the image of $\mathrm{ad}_{\Bar{\varphi}}$ and $S$ is indeed non-negative on the image of $\mathrm{ad}_{\Bar{\varphi}}$. The conclusion follows.
$\hfill\square$

\hspace*{\fill}

In the following proposition, we establish a lower bound on $\frac{\mathrm{d}}{\mathrm{d}t}(Ric(S)_{11}+Ric(S)_{22})$. The estimation given below is more precise than that in \cite{Brendle19}. 

The ideas for our improvement are as follows. In order to control $\frac{\mathrm{d}}{\mathrm{d}t}(Ric(S)_{11}+Ric(S)_{22})$, a crucial ingredient we need is to give an upper bound for $Ric(S)_{22}-Ric(S)_{11}$. By calculation, $\frac{\mathrm{d}}{\mathrm{d}t}(Ric(S)_{11}+Ric(S)_{22}+\frac{2\gamma}{n}\scal(S))$ contains a term in the form of $c|Ric_0(S)|^2$, where $c=c(b)$ is positive. One could use the direct estimate $|Ric_0(S)|^2\geq0$, but here we give a more precise estimate. Notice that if the quantity $Ric(S)_{22}-Ric(S)_{11}$ is large, then $|Ric_0(S)|^2$ cannot be too small, this inspires us to establish a lower bound of $|Ric_0(S)^2|$ in terms of $(Ric(S)_{22}-Ric(S)_{11})^2$, as we will show in the lemma \ref{lemma31}.

\begin{proposition}\label{prop 4th condition for first cone}
    Let $S\in\mathcal{E}(b)\backslash\{0\}$ and $Ric(S)_{11}+Ric(S)_{22}+\frac{2\gamma}{n}\scal(S)=0$ for some pair of orthonormal vectors $\{e_1,e_2\}$, then:
    \begin{equation}
        \frac{\mathrm{d}}{\mathrm{d}t}\left(Ric(S)_{11}+Ric(S)_{22}+\frac{2\gamma}{n}\scal(S)\right)>0.
    \end{equation}
\end{proposition}

\begin{lemma}\label{lemma31}
    For any $R\in S^2(\mathfrak{so}(n))$, let $Ric_{11}\leq Ric_{22}$ be the two smallest eigenvalues for $Ric(R)$. We have

\begin{equation}\label{elemma31}
    |Ric_0|^2\geq \frac{n-1}{n}(Ric_{22}-Ric_{11})^2.
\end{equation}
\end{lemma}

\noindent\textit{Proof.} We may assume $\lambda_1\leq\lambda_2\leq\cdots\leq\lambda_{n}$ are eigenvalues for $Ric$, and $\Bar{\lambda}=\frac{scal}{n}=\frac{1}{n}\sum_{i=1}^n\lambda_i$. Then we have
\begin{equation*}
    \begin{split}
    |Ric_0|^2&=\sum_{i=1}^n(\lambda_i-\Bar{\lambda})^2\\
    &\geq (\lambda_1-\Bar{\lambda})^2+(\lambda_2-\Bar{\lambda})^2+(n-2)\left(\frac{\sum_{i=3}^n\lambda_i}{n-2}-\Bar{\lambda}\right)^2.
    \end{split}
\end{equation*}

Now, if we replace $R$ by $R+c\mathrm{id}\wedge\mathrm{id}$ for some real number $c$, then the two sides of \ref{elemma31} remain the same. By the homogeneity of the inequality, we can assume without loss of generality that $\lambda_1=-1,\lambda_2=0$. Denote $\nu=\frac{1}{n-2}\sum_{i=3}^n\lambda_i$, then $\nu\geq0$ and $\Bar{\lambda}=\frac{1}{n}((n-2)\nu-1)$, we calculate that

\begin{equation*}
    \begin{split}
        |Ric_0|^2&\geq\frac{((n-2)\nu+(n-1))^2}{n^2}+\frac{((n-2)\nu-1)^2}{n^2}+(n-2)\frac{(2\nu+1)^2}{n^2}\\
        &=\frac{1}{n^2}(2n(n-2)\nu^2+2n(n-2)\nu+n(n-1))\\
        &\geq\frac{n-1}{n}.
    \end{split}
\end{equation*}

The assertion follows.$\hfill\square$

\hspace*{\fill}

\noindent\textit{Proof of Proposition \ref{prop 4th condition for first cone}.}  We may assume that $\{e_1,e_2\}$ are eigenvectors of $Ric(S)$, and extend $\{e_1,e_2\}$ to an eigenbasis $\{e_1,\cdots,e_n\}$ for $Ric(S)$. Let  $\lambda_{1}\leq\lambda_{2}\leq\cdots\leq\lambda_{n}$ denote the eigenvalues for $Ric(S)$ with respect to the basis $\{e_1,\cdots,e_n\}$. Following the calculations in \cite{Brendle19}, we have

\begin{equation}
    \begin{split}
        &\frac{1}{2}\frac{\mathrm{d}}{\mathrm{d}t}\left(Ric(S)_{11}+Ric(S)_{22}+\frac{2\gamma}{n}\scal(S)\right)\\
        =&\sum_{p=3}^n(S_{1p1p}+S_{2p2p})(\lambda_{p}-\frac{1}{2}\left(\lambda_2+\lambda_1)\right)\\
        &+\frac{4(1+\gamma)}{n^2}(a-b)\scal^2(S)+\frac{2\gamma(1+\gamma)}{n^2}(1-2b)(\scal(S))^2\\
        &-b(\lambda_2-\lambda_1)^2+\frac{2\gamma}{n}(1-2b)|Ric_0(S)|^2\\
        &+2(1+\gamma)\frac{n^2b^2-2(n-1)(a-b)(1-2b)}{n(1+2(n-1)a)}|Ric_0(S)|^2\\
    \geq&\sum_{p=3}^n(S_{1p1p}+S_{2p2p})(\lambda_{p}-\frac{1}{2}\left(\lambda_2+\lambda_1)\right)\\
    &+\frac{2(1+\gamma)}{n^2}\frac{b(1+(n-2)b)^2}{2+(n-3)b}(\scal(S))^2-\rho(\lambda_2-\lambda_1)^2.
    \end{split}
\end{equation}
where in the last inequality we used Lemma \ref{lemma31}, definition of $\rho$, and the identity
\begin{equation*}
    2(a-b)+\gamma(1-2b)=\frac{b(1+(n-2)b)^2}{2+(n-3)b}.
\end{equation*}

Observe that $S\in C_{PIC}$, we have
\begin{equation*}
    \scal(S)-2Ric(S)_{kk}=\sum_{\substack{i\neq k\\j\neq k}}S_{ijij}\geq0
\end{equation*}
for each $k\in\{1,2,\cdots,n\}$. It follows that each eigenvalue of $Ric(S)$ is bounded above from the sum of other eigenvalues. Combining with the fact that $\lambda_1+\lambda_2\leq 0$, this gives us the inequality

\begin{equation*}
    \sum_{p\in\{3,\cdots,n\}\backslash\{m\}}(\lambda_p-\frac{1}{2}(\lambda_2+\lambda_1))\geq\lambda_m-\frac{1}{2}(\lambda_2+\lambda_1).
\end{equation*}

As $S-T\in C_{PIC}$, it follows that there exists at most one $m\in\{3,\cdots,n\}$ such that $S_{1m1m}+S_{2m2m}-T_{1m1m}-T_{2m2m}<0$. Moreover, for any $p\in\{3,\cdots,n\}\backslash\{m\}$ we have $$(S_{1m1m}+S_{2m2m}-T_{1m1m}-T_{2m2m})+(S_{1p1p}+S_{2p2p}-T_{1p1p}-T_{2p2p})>0,$$so we have

\begin{equation*}
    \sum_{p=3}^n(S_{1p1p}=S_{2p2p}-T_{1p1p}-T_{2p2p})(\lambda_p-\frac{1}{2}(\lambda_1+\lambda_2))\geq0.
\end{equation*}

Since $T\geq0$, we have
\begin{equation*}
    \begin{split}
        &\sum_{p=3}^n(S_{1p1p}+S_{2p2p})(\lambda_p-\frac{1}{2}(\lambda_1+\lambda_2))\\
        \geq&\sum_{p=3}^n(T_{1p1p}+T_{2p2p})(\lambda_p-\frac{1}{2}(\lambda_1+\lambda_2))\\
        \geq&\frac{1}{2}\sum_{p=3}^n(T_{1p1p}+T_{2p2p})(\lambda_2-\lambda_1).
    \end{split}
\end{equation*}

It now follows that 
\begin{equation*}
    \begin{split}
        &\frac{1}{2}\frac{\mathrm{d}}{\mathrm{d}t}\left(Ric(S)_{11}+Ric(S)_{22}+\frac{2\gamma}{n}\scal(S)\right)\\
        \geq&\frac{1}{2}\sum_{p=3}^n(T_{1p1p}+T_{2p2p})(\lambda_2-\lambda_1)\\
        &+\frac{2(1+\gamma)}{n^2}\frac{b(1+(n-2)b)^2}{2+(n-3)b}\scal^2(S)-\rho(\lambda_2-\lambda_1)^2.
    \end{split}
\end{equation*}

Recall the fourth property in Definition \ref{def firstcone}, we have

\begin{equation*}
    \sum_{p=3}^n(T_{1p1p}+T_{2p2p})\geq \frac{(\lambda_2-\lambda_1)^2}{\omega\scal(S)}.
\end{equation*}

If $\lambda_2-\lambda_1=0$, then clearly we have $\frac{\mathrm{d}}{\mathrm{d}t}\left(Ric(S)_{11}+Ric(S)_{22}+\frac{2\gamma}{n}\scal(S)>0\right)$. Now assume $\lambda_2-\lambda_1>0$, the elementary inequality $2x+y\geq3(x^2y)^{\frac{1}{3}}$ for $x,y\geq0$ gives

\begin{equation*}
\frac{\lambda_2-\lambda_1}{2\omega\scal(S)}+\frac{2(1+\gamma)}{n^2}\frac{b(1+(n-2)b)^2}{2+(n-3)b}\frac{\scal^2(S)}{(\lambda_2-\lambda_1)^2}
        \geq3\left(\frac{1+\gamma}{8\omega^2n^2}\frac{b(1+(n-2)b)^2}{2+(n-3)b}\right)^{\frac{1}{3}}=(1-\varepsilon)^{-\frac{2}{3}}\rho>\rho,
\end{equation*}
where we have used the definition of $\gamma, \rho$ and $\omega$.

As a consequence,
\begin{equation*}
\begin{split}
    &\frac{1}{2}\frac{\mathrm{d}}{\mathrm{d}t}\left(Ric(S)_{11}+Ric(S)_{22}+\frac{2\gamma}{n}\scal(S)\right)\\
    \geq&\frac{(\lambda_2-\lambda_1)^3}{2\omega\scal(S)}+\frac{2(1+\gamma)}{n^2}\frac{b(1+(n-2)b)^2}{2+(n-3)b}(\scal(S))^2-\rho(\lambda_2-\lambda_1)^2\\
    >&0.
\end{split}
\end{equation*}

This proves the assertion.$\hfill\square$

\hspace*{\fill}

Finally, since we assumed an upper bound for $(Ric(S)_{22}-Ric(S)_{11})$ in Definition \ref{def firstcone}, we need to verify that this bound is preserved along the evolution $\frac{\mathrm{d}}{\mathrm{d}t}S=D_{a,b}(S)+Q(S)$ for $S$. The estimation below is more precise than that in \cite{Brendle19}. Our estimations are inspired by the expression formula of the evolution of scalar curvature
\begin{equation*}
    \frac{\mathrm{d}}{\mathrm{d}t}(\scal(S))=P|Ric(S)|^2+Q(\scal(S))^2,
\end{equation*}

as we shall see below.

\begin{proposition}\label{prop 5th condition of first cone}
    Suppose $S\in \mathcal{E}(b)\backslash\{0\}$, $T\geq0$ and $S-T\in C_{PIC}$. Moreover, suppose that
    \begin{equation*}
            Ric(S)_{22}-Ric(S)_{11}=\omega^{\frac{1}{2}}(\scal(S))^{\frac{1}{2}}\left(\sum_{p=3}^n(T_{1p1p}+T_{2p2p})\right)^{\frac{1}{2}},
    \end{equation*}

    then
    \begin{equation}
        \frac{\mathrm{d}}{\mathrm{d}t}(Ric(S)_{22}-Ric(S)_{11})<\frac{\mathrm{d}}{\mathrm{d}t}\left[\omega^{\frac{1}{2}}(\scal(S))^{\frac{1}{2}}\left(\sum_{p=3}^n(T_{1p1p}+T_{2p2p})\right)^{\frac{1}{2}}\right]
    \end{equation}
\end{proposition}

\noindent\textit{Proof}. 
    We assume $\sum_{p=3}^n(T_{1p1p}+T_{2p2p})>0$ (otherwise the right hand side of the desired inequality is infinity and the result is trivially true), clearly we have $Ric(S)_{22}-Ric(S)_{11}\geq0$ and $Ric(S)_{12}=0$. Without loss of generality, we may assume that $Ric(S)_{pq}=0$ for $3\leq p<q\leq n$. As in \cite{Brendle19} we write
    \begin{equation*}
        \frac{\mathrm{d}}{\mathrm{d}t}(Ric(S)_{22}-Ric(S)_{11})=J_{1}+J_{2}+J_{3}+J_{4}+J_{5}+J_{6},
    \end{equation*}
where the terms $J_{1},\cdots, J_{6}$ are defined by
\begin{equation*}
    \begin{split}
        J_{1}&=2\sum_{p=3}^n(S_{2p2p}-S_{1p1p})Ric(S)_{pp},\\
        J_{2}&=-4\sum_{p=3}^nS_{121p}Ric(S)_{2p},\\
        J_{3}&=4\sum_{p=3}^nS_{212p}Ric(S)_{1p},\\
        J_{4}&=-2S_{1212}(Ric(S)_{22}-Ric(S)_{11}),\\
        J_{5}&=4b((Ric(S)^2)_{11}-(Ric(S)^2)_{22}),\\
        J_{6}&=\frac{4}{n}(2b+(n-2)a)\scal(S)(Ric(S)_{22}-Ric(S)_{11}).\\
    \end{split}
\end{equation*}

It turns out that the terms $J_1$ and $J_6$ play a major role, and it's more convenient to put the terms $J_4, J_5, J_6$ together.

We estimate the terms $J_{1},J_{2},J_{3}$ by 
\begin{equation*}
\begin{split}
    J_{1}&\leq \tau\sum_{p=3}^n(S_{2p2p}-S_{1p1p})^2+\frac{1}{\tau}\sum_{p=3}^n(Ric(S)_{pp})^2,\\
    J_{2}&\leq 2\tau(1+b\sqrt{n-2})\sum_{p=3}^n(S_{121p})^2+\frac{2}{\tau(1+b\sqrt{n-2})}\sum_{p=3}^n(Ric_{2p})^2,\\
    J_{3}&\leq 2\tau(1+b\sqrt{n-2})\sum_{p=3}^n(S_{212p})^2+\frac{2}{\tau(1+b\sqrt{n-2})}\sum_{p=3}^n(Ric_{1p})^2.
\end{split}
\end{equation*}

Notice that
\begin{equation*}
    \sum_{3\leq p\neq q\leq n}(S_{1p1p}+S_{1q1q}+S_{2p2p}+S_{2q2q})\geq 0,
\end{equation*}
we obtain $Ric(S)_{11}+Ric(S)_{22}\geq 2S_{1212}$. Combining this with the assumption that $Ric(S)_{22}-Ric(S)_{11}\geq0$, the term $J_{5}$ can be estimated by

\begin{equation*}
    \begin{split}
        J_5\leq& -4b(Ric(S)_{22}-Ric(S)_{11})(Ric(S)_{11}+Ric(S)_{22})+4b\sum_{p=3}^n(Ric(S)_{1p})^2\\
        \leq& -8b(Ric(S)_{22}-Ric(S)_{11})S_{1212}+2b\tau\left(\frac{1}{\sqrt{n-2}}+b\right)\sum_{p=3}^n(Ric(S)_{1p})^2\\
        &+\frac{2b}{\tau\left(\frac{1}{\sqrt{n-2}}+b\right)}\sum_{p=3}^n(Ric(S)_{1p})^2\\
        =& -8b(Ric(S)_{22}-Ric(S)_{11})S_{1212}+2b\tau\left(\frac{1}{\sqrt{n-2}}+b\right)\sum_{p=3}^n\left(S_{12p2}+\sum_{q=3}^nS_{1qpq}\right)^2\\
        &+\frac{2b}{\tau\left(\frac{1}{\sqrt{n-2}}+b\right)}\sum_{p=3}^n(Ric(S)_{1p})^2\\
        \leq& -8b(Ric(S)_{22}-Ric(S)_{11})S_{1212}+2(n-2)b\tau\left(\frac{1}{\sqrt{n-2}}+b\right)\sum_{p=3}^n\left((S_{12p2})^2+\sum_{q=3}^n(S_{1qpq})^2\right)\\
        &+\frac{2b}{\tau\left(\frac{1}{\sqrt{n-2}}+b\right)}\sum_{p=3}^n(Ric(S)_{1p})^2\\
        :=&-8bS_{1212}(Ric(S)_{22}-Ric(S)_{11})+J_{7}.
    \end{split}
\end{equation*}

Notice that we have the inequality

\begin{equation*}
        \begin{split}
            0\leq&\sum_{3\leq p\neq q\leq n}(S_{1212}+S_{1p1p}+S_{2q2q}+S_{pqpq})\\
            &=(n-2)(n-3)S_{1212}+(n-4)\sum_{p=3}^n(S_{1p1p}+S_{2p2p})+\sum_{p=3}^n Ric(S)_{pp}\\
            &=(n-1)(n-4)S_{1212}+(n-5)\sum_{p=3}^n(S_{1p1p}+S_{2p2p})+\scal(S).\\
        \end{split}
\end{equation*}

This allows us to estimate $J_{4}+J_{5}+J_{6}$, as follows:

\begin{equation*}
    \begin{split}
        J_4+&J_5+J_6\\
        \leq&\left((-2-8b)S_{1212}+\frac{4}{n}(2b+(n-2)a)\scal(S)\right)(Ric(S)_{22}-Ric(S)_{11})+J_7\\
        \leq&\left[\frac{(2+8b)(n-5)}{(n-1)(n-4)}\sum_{p=3}^n(S_{1p1p}+S_{2p2p})+\left(\frac{2+8b}{(n-1)(n-4)}+\frac{4}{n}(2b+(n-2)a)\right)\scal(S)\right]\\
        & \cdot (Ric(S)_{22}-Ric(S)_{11})+J_7\\
        =&\frac{(2+8b)(n-5)}{(n-1)(n-4)}\left(\sum_{p=3}^n(S_{1p1p}+S_{2p2p})\right)(Ric(S)_{22}-Ric(S)_{11})\\
        &+A\scal(S)\left(\sum_{p=3}^n(S_{1p1p}-S_{2p2p})\right)+J_7,\\
    \end{split}
\end{equation*}
where we recall that in the beginning of this section, we define $A=\frac{2+8b}{(n-1)(n-4)}+\frac{4}{n}(2b+(n-2)a)$. Thus

\begin{equation*}
    \begin{split}
    J_4+&J_5+J_6\\
        \leq & \left(1+\frac{A}{2}\sqrt{n(n-2)}\right)\tau\sum_{p=3}^n(S_{1p1p}+S_{2p2p})^2+\frac{C}{\tau}(Ric(S)_{22}-Ric(S)_{11})^2\\
        &+\frac{A\tau}{2}\sqrt{n(n-2)}\sum_{p=3}^n(S_{2p2p}-S_{1p1p})^2+\frac{A}{2\tau}\sqrt{\frac{n-2}{n}}(\scal(S))^2+J_7,\\
    \end{split}
\end{equation*}
where
\begin{equation*}
    C=\frac{(1+4b)^2(n-5)^2(n-2)}{(n-1)^2(n-4)^2\left(1+\frac{A}{2}\sqrt{n(n-2)}\right)}\leq\frac{1+4b}{n-1}<\frac{1}{2}.
\end{equation*}

Summing up, we have 
\begin{equation*}
    \begin{split}
        \frac{\mathrm{d}}{\mathrm{d}t}&(Ric(S)_{22}-Ric(S)_{11})\\
        \leq& \left(1+\frac{A}{2}\sqrt{n(n-2)}\right)\tau\sum_{p=3}^n(S_{2p2p}-S_{1p1p})^2+\left(1+\frac{A}{2}\sqrt{n(n-2)}\right)\tau\sum_{p=3}^n(S_{1p1p}+S_{2p2p})^2\\
        &+2\tau(1+b\sqrt{n-2})\sum_{p=3}^n(S_{121p})^2+\left(2\tau(1+b\sqrt{n-2})+2(n-2)b\tau\left(\frac{1}{\sqrt{n-2}}+b\right)\right)\sum_{p=3}^n(S_{212p})^2\\
        &+2(n-2)b\tau\left(\frac{1}{\sqrt{n-2}}+b\right)\sum_{p=3}^n\sum_{q=3}^n(S_{1qpq})^2\\
        &+\frac{1}{\tau}\sum_{p=3}^n(Ric(S)_{pp})^2+\frac{A}{2\tau}\sqrt{\frac{n-2}{n}}(\scal(S))^2+\frac{1}{2\tau}(Ric(S)_{22}-Ric(S)_{11})^2\\
        &+\frac{2}{\tau(1+b\sqrt{n-2})}\sum_{p=3}^n(Ric(S)_{2p})^2+\left(\frac{2}{\tau(1+b\sqrt{n-2})}+\frac{2b}{\tau\left(\frac{1}{\sqrt{n-2}}+b\right)}\right)\sum_{p=3}^n(Ric(S)_{1p})^2.
    \end{split}
\end{equation*}

Notice that 

\begin{equation*}
    \frac{1}{1+b\sqrt{n-2}}+\frac{b}{\frac{1}{\sqrt{n-2}}+b}=1,
\end{equation*}
and by Lemma \ref{lem computational results for first family}, we have
\begin{equation*}
    2(1+b\sqrt{n-2})+2(n-2)b\left(\frac{1}{\sqrt{n-2}}+b\right)=2(1+b\sqrt{n-2})^2\leq 2+A\sqrt{n(n-2)}.
\end{equation*}

It follows that

\begin{equation*}
    \begin{split}
        \frac{\mathrm{d}}{\mathrm{d}t}&(Ric(S)_{22}-Ric(S)_{11})\\
       \leq &(2+A\sqrt{n(n-2)})\tau\sum_{p=3}^n\left((S_{1p1p})^2+(S_{2p2p})^2+(S_{121p})^2+(S_{212p})^2+\sum_{q=3}^n(S_{1qpq})^2\right)\\
       &+\frac{1}{\tau}\sum_{p=3}^n((Ric(S)_{pp})^2+2(Ric(S)_{1p})^2+2(Ric(S)_{2p})^2)+\frac{1}{\tau}((Ric(S)_{11})^2+(Ric(S)_{22})^2)\\
       &+\frac{A}{2\tau}\sqrt{\frac{n-2}{n}}(\scal(S))^2\\
       \leq & \left(1+\frac{A}{2}\sqrt{n(n-2)}\right)\tau\sum_{p=3}^n((S^2)_{1p1p}+(S^2)_{2p2p})+\frac{1}{\tau}|Ric(S)|^2+\frac{A}{2\tau}\sqrt{\frac{n-2}{n}}(\scal(S))^2,
    \end{split}
\end{equation*}
for any $\tau>0$.

At this point, we make the following claim:

\noindent\textbf{Claim.} We have for any $\sigma>0$ the estimate:
\begin{equation*}
    \frac{\mathrm{d}}{\mathrm{d}t}(Ric(S)_{22}-Ric(S)_{11})<\frac{\omega^{\frac{1}{2}}\sigma}{2}\sum_{p=3}^n((S^2)_{1p1p}+(S^2)_{2p2p})+\frac{\omega^{\frac{1}{2}}}{2\sigma}(P|Ric(S)|^2+Q(\scal(S))^2).
\end{equation*}

We shall distinguish 2 cases:

\noindent\textit{Case 1}: If $\frac{AP}{2}\sqrt{\frac{n-2}{n}}-Q\geq0$, we split the term containing $(\scal(S))^2$ into two parts:

\begin{equation*}
    \frac{A}{2\tau}\sqrt{\frac{n-2}{n}}=\left(\frac{AP}{2(P+nQ)\tau}\sqrt{\frac{n-2}{n}}-\frac{Q}{(P+nQ)\tau}\right)+\left(\frac{AQ\sqrt{n(n-2)}}{2(P+nQ)\tau}+\frac{Q}{(P+nQ)\tau}\right).
\end{equation*}

Applying the inequality $(\scal(S))^2\leq n|Ric(S)|^2$, we obtain
\begin{equation*}
    \begin{split}
        \frac{\mathrm{d}}{\mathrm{d}t}&(Ric(S)_{22}-Ric(S)_{11})\\
        \leq&\left(1+\frac{A}{2}\sqrt{n(n-2)}\right)\tau\sum_{p=3}^n((S^2)_{1p1p}+(S^2)_{2p2p})\\
        &+\left(\frac{AP\sqrt{n(n-2)}}{2(P+nQ)\tau}-\frac{nQ}{(P+nQ)\tau}+\frac{1}{\tau}\right)|Ric(S)|^2\\
        &+\left(\frac{AQ\sqrt{n(n-2)}}{2(P+nQ)\tau}+\frac{Q}{(P+nQ)\tau}\right)(\scal(S))^2\\
        =&\left(1+\frac{A}{2}\sqrt{n(n-2)}\right)\tau\sum_{p=3}^n((S^2)_{1p1p}+(S^2)_{2p2p})\\
        &+\left(1+\frac{A}{2}\sqrt{n(n-2)}\right)\frac{1}{(P+nQ)\tau}(P|Ric(S)|^2+Q(\scal(S))^2).
    \end{split}
\end{equation*}

Recall that Lemma \ref{lem computational results for first family} shows that
\begin{equation*}
    \left(1+\frac{A}{2}\sqrt{n(n-2)}\right)^2\frac{1}{P+nQ}<\frac{\omega}{4}.
\end{equation*}

Thus, we have
\begin{equation*}
    \frac{\mathrm{d}}{\mathrm{d}t}(Ric(S)_{22}-Ric(S)_{11})<\frac{\omega^{\frac{1}{2}}\sigma}{2}\sum_{p=3}^n((S^2)_{1p1p}+(S^2)_{2p2p})+\frac{\omega^{\frac{1}{2}}}{2\sigma}(P|Ric(S)|^2+Q(\scal(S))^2),
\end{equation*}
for any $\sigma>0$, the claim is true.

\noindent\textit{Case 2}: If otherwise $\frac{AP}{2}\sqrt{\frac{n-2}{n}}-Q<0$, then

\begin{equation*}
    \begin{split}
        \frac{\mathrm{d}}{\mathrm{d}t}&(Ric(S)_{22}-Ric(S)_{11})\\
        \leq&\left(1+\frac{A}{2}\sqrt{n(n-2)}\right)\tau\sum_{p=3}^n((S^2)_{1p1p}+(S^2)_{2p2p})+\frac{1}{\tau}|Ric(S)|^2+\frac{Q}{P\tau}(\scal(S))^2\\
        &=\left(1+\frac{A}{2}\sqrt{n(n-2)}\right)\tau\sum_{p=3}^n((S^2)_{1p1p}+(S^2)_{2p2p})+\frac{1}{P\tau}(P|Ric(S)|^2+Q(\scal(S))^2).
    \end{split}
\end{equation*}

Recall that Lemma \ref{lem computational results for first family} shows that
\begin{equation*}
    \left(1+\frac{A}{2}\sqrt{n(n-2)}\right)\frac{1}{P}<\frac{\omega}{4},
\end{equation*}
thus for any $\sigma>0$, we have
\begin{equation*}
    \frac{\mathrm{d}}{\mathrm{d}t}(Ric(S)_{22}-Ric(S)_{11})<\frac{\omega^{\frac{1}{2}}\sigma}{2}\sum_{p=3}^n((S^2)_{1p1p}+(S^2)_{2p2p})+\frac{\omega^{\frac{1}{2}}}{2\sigma}(P|Ric(S)|^2+Q(\scal(S))^2),
\end{equation*}
also proves the claim.

Arriving here, combining the fact that$\frac{\mathrm{d}}{\mathrm{d}t}(\scal(S))=P|Ric(S)|^2+Q(\scal(S))^2$ with the inequality
\begin{equation*}
    \frac{\mathrm{d}}{\mathrm{d}t}\sum_{p=3}^n(T_{1p1p}+T_{2p2p})>\sum_{p=3}^n((S^2)_{1p1p}+(S^2)_{2p2p}),
\end{equation*}
we are able to deduce that

    \begin{align*}
        \frac{\mathrm{d}}{\mathrm{d}t}&\left[\omega^{\frac{1}{2}}\scal^{\frac{1}{2}}\left(\sum_{p=3}^n(T_{1p1p}+T_{2p2p})\right)^{\frac{1}{2}}\right] \\
        =&\frac{1}{2}\omega^{\frac{1}{2}}\scal(S)^{\frac{1}{2}}\left(\sum_{p=3}^n(T_{1p1p}+T_{2p2p})\right)^{-\frac{1}{2}}\frac{\mathrm{d}}{\mathrm{d}t}\left(\sum_{p=3}^n(T_{1p1p}+T_{2p2p})\right)\\
        &+\frac{1}{2}\omega^{\frac{1}{2}}\scal(S)^{-\frac{1}{2}}\left(\sum_{p=3}^n(T_{1p1p}+T_{2p2p})\right)^{\frac{1}{2}}\frac{\mathrm{d}}{\mathrm{d}t}\scal(S)\\
        >&\frac{\omega^{\frac{1}{2}}\sigma}{2}\sum_{p=3}^n((S^2)_{1p1p}+(S^2)_{2p2p})+\frac{\omega^{\frac{1}{2}}}{2\sigma}(P|Ric(S)|^2+Q(\scal(S))^2),\\
    \end{align*}
where we take $\sigma=\scal(S)^{\frac{1}{2}}\left(\sum_{p=3}^n(T_{1p1p}+T_{2p2p})\right)^{-\frac{1}{2}}$. Putting these facts together, the assertion follows.$\hfill\square$

\hspace*{\fill}

We see that Proposition \ref{prop 2ndprop of first cone} to Proposition \ref{prop 5th condition of first cone} guarantees the validity of Theorem \ref{thm goal for first family}.

\section{Proof of Theorem 1.2}

In this section, we construct another one-parameter family of cones $\{\Tilde{C}(b)\},0
<b\leq \Tilde{b}_{max}$ that are invariant under Hamiltonian ODE. Here we adapt the ideas of construction in Section 4 of $\cite{Brendle19}$ and change some related numerical data. This family $\Tilde{C}(b)$ pinches towards $C_{PIC1}$ and can be joined continuous to the family $C(b)$ constructed in the previous section. We concatenate the two families to get the desired family of cones required in Theorem \ref{thm main cone}.

We assume in this section that $9\leq n\leq 11$. We define $\Tilde{b}_{max}=\frac{1}{5n}$. The following definition is from Definition 4.1 of \cite{Brendle19}:

\begin{definition}\label{def second family}
    Assume $0<b\leq \Tilde{b}_{max}$ and let $a=b+\frac{n-2}{2}b^2$. We denote by $\Tilde{\mathcal{E}}(b)$ the set of all $S\in S_{B}^2(\mathfrak{so}(n))$ such that $l_{a,b}(S)\in C(b_{max})$ and
    \begin{equation*}
        Z:=S_{1313}+\lambda^2S_{1414}+S_{2323}-2\lambda S_{1234}+\sqrt{2a}(1-\lambda^2)(Ric(S)_{11}+Ric(S)_{22})\geq 0
    \end{equation*}
    for every otrhonormal four-frame $\{e_1,e_2,e_3,e_4\}$ and every $\lambda\in[0,1]$. We then define $\Tilde{C}(b)=l_{a,b}(\Tilde{\mathcal{E}}(b))$.
\end{definition}

It's clear that $\Tilde{C}(b)$ is convex for each $b$. We now verify that the family $C(b), 0<b\leq b_{max}$, constructed in the previous section can be joined with the family $\Tilde{C}(b), 0<b\leq \Tilde{b}_{max}$. Recall that $b_{max}=\frac{1}{2n+2}$. We define

\begin{equation*}
    a_{max}=\frac{(2+(n-2)b_{max})^2}{2(2+(n-3)b_{max})}b_{max}\quad and \quad \gamma_{max}=\frac{b_{max}}{2+(n-3)b_{max}}.
\end{equation*}

\subsection{Elementary lemmas}

By direct calculation, we can verify the following two elementary lemmas are true.

\begin{lemma}\label{lemma connect cone inequalities}(cf. Lemma 4.2 of \cite{Brendle19})
    Let $\Tilde{b}_{max}$ be defined as above, and let $\Tilde{a}_{max}$ be defined by $\Tilde{a}_{max}=\Tilde{b}_{max}+\frac{n-2}{2}\Tilde{b}_{max}^2$. Then

    \begin{equation}\label{equ 1 to connect cone}
        \frac{1+(n-2)b_{max}}{1+(n-2)\Tilde{b}_{max}}\sqrt{2\Tilde{a}_{max}}\geq\frac{n^2-5n+4}{n^2-7n+14}\cdot\frac{1}{n-4},
    \end{equation}

    \begin{equation}\label{equ 2 to connect cone}
        \frac{a_{max}-\Tilde{a}_{max}}{1+2(n-1)\Tilde{a}_{max}}-\frac{b_{max}-\Tilde{b}_{max}}{1+(n-2)\Tilde{b}_{max}}\geq 0,
    \end{equation}
and
    \begin{equation}\label{equ 3 to connect cone}
        \begin{split}
            &2\left(\frac{a_{max}-\Tilde{a}_{max}}{1+2(n-1)\Tilde{a}_{max}}-(1+\gamma_{max})\frac{b_{max}-\Tilde{b}_{max}}{1+(n-2)\Tilde{b}_{max}}\right)\\
            &+\left(\frac{2(n-1)(a_{max}-\Tilde{a}_{max})}{1+2(n-1)\Tilde{a}_{max}}-\frac{(n-2)(b_{max}-\Tilde{b}_{max})}{1+(n-2)\Tilde{b}_{max}}\right)\sqrt{2\Tilde{a}_{max}}\\
            \geq&\frac{1+(n-2)b_{max}}{1+(n-2)\Tilde{b}_{max}}\cdot\frac{n\sqrt{2\Tilde{a}_{max}}}{n^2-5n+4}.
        \end{split}
    \end{equation}
\end{lemma}

\textit{Proof.} We verify by direct calculations the above formulas for $n=9,10,11$ respectively.
For the equation \ref{equ 1 to connect cone}, we have
    \begin{itemize}
        \item For $n=9$, $LHS\approx 0.255692>\frac{1}{4} =RHS$;
        \item  For $n=10$, $LHS\approx 0.244333>\frac{9}{44} =RHS$;
        \item For $n=11$, $LHS\approx 0.234367>\frac{5}{29} =RHS$.
    \end{itemize}

For the equation \ref{equ 2 to connect cone}, we have
    \begin{itemize}
        \item For $n=9$, $LHS\approx 0.002043>0$;
        \item  For $n=10$, $LHS\approx 0.001942>0$;
        \item For $n=11$, $LHS\approx 0.001846>0$.
    \end{itemize}

For the equation \ref{equ 3 to connect cone}, we have
    \begin{itemize}
        \item For $n=9$, $LHS\approx 0.057547>0.057531 \approx RHS$;
        \item  For $n=10$, $LHS\approx 0.055897>0.045246 \approx RHS$;
        \item For $n=11$, $LHS\approx 0.054324>0.036829 \approx RHS$.
    \end{itemize}
$\hfill\square$

\begin{lemma}\label{lemma elementary second cone invariant}(cf. Lemma 4.4 of \cite{Brendle19})
    Assume that $0<b\leq \Tilde{b}_{max}$, and let $a$ be defined by $a=b+\frac{n-2}{2}b^2$. Then
    \begin{equation*}
        nb^2(1-2b)-2(a-b)(1-2b+nb^2)\geq0,
    \end{equation*}
    \begin{equation*}
        n^2b^2-2(n-1)(a-b)(1-2b)\geq0.
    \end{equation*}

Moreover, if we put
\begin{equation*}
    \zeta:=\frac{1+2(n-1)a}{1+2(n-1)a_{max}}\cdot\frac{1+(n-2)b_{max}}{1+(n-2)b}(1+\gamma_{max}),
\end{equation*}
then $\zeta\leq1$ and $1+(n-2)(1-\zeta)\geq2\zeta^2\frac{n^2-2n+2}{(n-2)^2}$.
\end{lemma}

\noindent\textit{Proof.} 
The first and second statements are direct since

\begin{equation*}
    nb^2(1-2b)-2(a-b)(1-2b+nb^2)=b^2(2-4b-n(n-2)b^2)\geq0,
\end{equation*}
and
\begin{equation*}
    n^2b^2--2(n-1)(a-b)(1-2b)=b^2(3n-2+2(n-1)(n-2)b)\geq0.
\end{equation*}

For the third and fourth statements, write $\Tilde{a}_{max}=\Tilde{b}_{max}+\frac{n-2}{2}\Tilde{b}_{max}$. Then $\frac{1+2(n-1)a}{1+(n-2)b}\leq \frac{1+2(n-1)\Tilde{a}_{max}}{1+(n-2)\Tilde{b}_{max}}$, hence

\begin{equation*}
    \zeta\leq\zeta_{max}:=\frac{1+2(n-1)\Tilde{a}_{max}}{1+2(n-1)a_{max}}\cdot\frac{1+(n-2)b_{max}}{1+(n-2)\Tilde{b}_{max}}(1+\gamma_{max}).
\end{equation*}

Now we have $\zeta_{max}\leq 1$ and $1+(n-2)(1-\zeta_{max})\geq2\zeta_{max}^2\frac{n^2-2n+2}{(n-2)^2}$ for $9\leq n\leq 11$, since

\begin{itemize}
    \item For $n=9$, $\zeta_{max}\approx0.824287<1,$ and $1+(n-2)(1-\zeta_{max})-2\zeta_{max}^2\frac{n^2-2n+2}{(n-2)^2}\approx0.427367>0$;
    \item For $n=10$, $\zeta_{max}\approx0.822096<1,$ and $1+(n-2)(1-\zeta_{max})-2\zeta_{max}^2\frac{n^2-2n+2}{(n-2)^2}\approx0.691388>0$;
    \item For $n=11$, $\zeta_{max}\approx0.820267<1,$ and $1+(n-2)(1-\zeta_{max})-2\zeta_{max}^2\frac{n^2-2n+2}{(n-2)^2}\approx0.939657>0$;
\end{itemize}

It follows that $\zeta\leq 1$ and $1+(n-2)(1-\zeta)\geq2\zeta^2\frac{n^2-2n+2}{(n-2)^2}$. $\hfill\square$

\hspace*{\fill}

The data $\zeta$ gives a characterization of a lower bound on the sum of two smallest eigenvalues of $Ric(S)$ for $S\in\Tilde{\mathcal{E}}(b)$, as follows:

\begin{lemma}\label{lemma elemenary2 second cone invariant}(cf. Lemma 4.5 of \cite{Brendle19})
    Suppose $S\in\Tilde{\mathcal{E}}(b)$ for $0<b\leq \Tilde{b}_{max}$ and 
    \begin{equation*}
        \zeta:=\frac{1+2(n-1)a}{1+2(n-1)a_{max}}\frac{1+(n-2)b_{max}}{1+(n-2)b}(1+\gamma_{max})\leq 1.
    \end{equation*}

    Then we have
    \begin{equation*}
        Ric(S)_{11}+Ric(S)_{22}\geq\frac{2(1-\zeta)}{n}\scal(S).
    \end{equation*}
\end{lemma}

\noindent\textit{Proof.} Let $T=l_{a_{max},b_{max}}^{-1}(l_{a,b}(S))\in\mathcal{E}(b_{max})$, then use the property of the map $l_{a,b}$, we have $RiC_0(S)=\frac{1+(n-2)b_{max}}{1+(n-2)b}Ric_0(T)$ and $\scal(S)=\frac{1+2(n-1)a_{max}}{1+2(n-1)a}\scal(T)$. Recall Definition \ref{def firstcone} that $Ric_0(T)_{11}+Ric_0(T)_{22}+\frac{2(1+\gamma_{max})}{n}\scal(T)\geq0$, we obtain that 

\begin{equation*}
    Ric_0(S)_{11}+Ric_0(S)_{22}+\frac{2\zeta}{n}\scal(S)\geq0,
\end{equation*}

the conclusion follows. $\hfill\square$

\subsection{Proof of Theorem \ref{thm main cone}}

With the above elementary lemmas, we now establish that the two family of pinching cones can be glued together, and that for each $0<b\leq \Tilde{b}_{max}$, the cone $\Tilde{C}(b)$ is transversally invariant under the Hamiltonian ODE. If these two statements are true, then we can connect the family $\{C(b)\},0<b\leq b_{max}$ and $\{\Tilde{C}(b)\}, 0<b\leq \Tilde{b}_{max}$ together to obtain the desired family of pinching cones required in Theorem \ref{thm main cone}. These two statements are established in Proposition \ref{prop for cone connection} and Theorem \ref{thm charactrization of second cone}. The proofs here are analogous to that in \cite{Brendle19}. In fact, only the validity of Lemma \ref{lemma connect cone inequalities} and Lemma \ref{lemma elementary second cone invariant} requires a constraint for the dimension $n$. Other parts of the proof is valid for all $n\geq4$. We present the proof below for completeness. 

\begin{proposition}\label{prop for cone connection}(cf. Proposition 4.3 of \cite{Brendle19})
    Let $b_{max}=\frac{1}{2n+2}$ and $\Tilde{b}_{max}=\frac{1}{5n}$ be defined as above, then $\Tilde{C}(\Tilde{b}_{max})=C(b_{max})$.
\end{proposition}

\textit{Proof.} We abbreviate that $b=\Tilde{b}_{max}$ and $a=\Tilde{a}_{max}$. It is clear by definition that $\Tilde{C}(b)\subseteq C(b_{max})$, it suffices to show the other side of inclusion. To this end, we need to show that $l_{a,b}^{-1}(C(b_{max}))\in\Tilde{\mathcal{E}}(b)$. That is, for any algebraic curvature operator $S$ satisfying $l_{a,b}(S)\in C(b_{max})$, we need to show that
\begin{equation*}
    Z:=S_{1313}+\lambda^2S_{1414}+S_{2323}+\lambda^2S_{2424}-2\lambda S_{1234}+\sqrt{2a}(1-\lambda^2)(Ric(S)_{11}+Ric(S)_{22})\geq0
\end{equation*}
for every orthonormal four-frame $\{e_1,e_2,e_3,e_4\}$ and every $\lambda\in[0,1]$.  

Let $T=l_{a_{max},b_{max}}^{-1}(l_{a,b}(S))\in \mathcal{E}(b_{max})$, then by Definition \ref{def firstcone}, $T\in C_{PIC}$ and $Ric(T)_{11}+Ric(T)_{22}+\frac{2\gamma_{max}}{n}\scal(T)\geq0$. Then we have

\begin{equation*}
        S=T+\frac{b_{max}-b}{1+(n-2)b}Ric(T)\wedge\mathrm{id}+\frac{1}{n}\left(\frac{a_{max}-a}{1+2(n-1)a}-\frac{b_{max}-b}{1+(n-2)b}\right)\scal(T)\mathrm{id}\wedge\mathrm{id}\\
\end{equation*}

Hence

\begin{equation*}
    \begin{split}
        Z=&T_{1313}+\lambda^2T_{1414}+T_{2323}+\lambda^2T_{2424}-2\lambda T_{1234}\\
        &+\frac{b_{max}-b}{1+(n-2)b}((1+\lambda^2)Ric(T)_{11}+(1+\lambda^2)Ric(T)_{22}+2Ric(T)_{33}+2\lambda^2Ric(T)_{44})\\
        &+\frac{4}{n}\left(\frac{a_{max}-a}{1+2(n-1)a}-\frac{b_{max}-b}{1+(n-22)b}\right)(1+\lambda^2)\scal(T)\\
        &+\frac{1+(n-2)b_{max}}{1+(n-2)b}\sqrt{2a}(1-\lambda^2)(Ric(T)_{11}+Ric(T)_{22})\\
        &+\frac{2}{n}\left(\frac{2(n-1)(a_{max}-a)}{1+2(n-11)a}-\frac{(n-2)(b_{max}-b)}{1+(n-2)b}\right)\sqrt{2a}(1-\lambda^2)\scal(T).
    \end{split}
\end{equation*}

Now, as $T\in C_{PIC}$, consider the orthonormal four frame $\{e_1,e_2,e_3,\lambda e_{4}+\sqrt{1-\lambda^2}e_{p}\}$ for each $5\leq p\leq n$, we obtain

\begin{equation*}
    T_{1313}+\lambda^2T_{1414}+T_{2323}+\lambda^2T_{2424}-2\lambda T_{1234}+(1-\lambda^2)(T_{1p1p}+T_{2p2p})\geq0.
\end{equation*}

Summing up the above equation over $5\leq p\leq n$ and notice that

\begin{equation}
    \sum_{p=5}^n(T_{1p1p}+T_{2p2p})\leq\sum_{p=3}^n(T_{1p1p}+T_{2p2p})=Ric(T)_{11}+Ric(T)_{22}-2T_{1313},
\end{equation}
we obtain

\begin{equation*}
    (n-4)(T_{1313}+\lambda^2T_{1414}+T_{2323}+\lambda^2T_{2424}-2\lambda T_{1234})+(1-\lambda^2)(Ric(T)_{11}+Ric(T)_{22}-2T_{1212})\geq0.
\end{equation*}

Consequently, we have

\begin{equation*}
    \begin{split}
        Z\geq& -\frac{1-\lambda^2}{n-4}(Ric(T)_{11}+Ric(T)_{22}-2T_{1212})\\
        &+\frac{b_{max}-b}{1+(n-2)b}((1+\lambda^2)Ric(T)_{11}+(1+\lambda^2)Ric(T)_{22}+2Ric(T)_{33}+2\lambda^2Ric(T)_{44})\\
        &+\frac{4}{n}\left(\frac{a_{max}-a}{1+2(n-1)a}-\frac{b_{max}-b}{1+(n-22)b}\right)(1+\lambda^2)\scal(T)\\
        &+\frac{1+(n-2)b_{max}}{1+(n-2)b}\sqrt{2a}(1-\lambda^2)(Ric(T)_{11}+Ric(T)_{22})\\
        &+\frac{2}{n}\left(\frac{2(n-1)(a_{max}-a)}{1+2(n-11)a}-\frac{(n-2)(b_{max}-b)}{1+(n-2)b}\right)\sqrt{2a}(1-\lambda^2)\scal(T).\\
        =:&RHS.
    \end{split}
\end{equation*}

We now show that $RHS\geq0$ for all $\lambda\in [0,1]$. Since $RHS$ is a linear function of $\lambda^2$, we only need to verify the cases where $\lambda=1$ and $\lambda=0$.

\noindent\textit{Case 1.} Suppose first that $\lambda=1$, then since $Ric(T)_{11}+Ric(T)_{22}+Ric(T)_{33}+Ric(T)_{44}\geq0$ and $\scal(T)\geq0$, we have
\begin{equation*}
    \begin{split}
        RHS=&2\frac{b_{max}-b}{1+(n-2)b}(Ric(T)_{11}+Ric(T)_{22}+Ric(T)_{33}+Ric(T)_{44})\\
        &+\frac{8}{n}\left(\frac{a_{max}-a}{1+2(n-1)a}-\frac{b_{max}-b}{1+(n-2)b}\scal(T)\right)\\
        \geq0.
    \end{split}
\end{equation*}

\textit{Case 2.} Suppose then that $\lambda=0$, by Lemma \ref{lemma connect cone inequalities}, we have
\begin{equation*}
    \frac{1}{n-4}\leq \frac{1+(n-2)b_{max}}{1+(n-2)b}\sqrt{2a}\frac{n^2-7n+14}{n^2-5n+4}.
\end{equation*}

Using inequalities $Ric(T)_{11}+Ric(T)_{22}-2T_{1212}\geq0$ and 
\begin{equation*}
    Ric(T)_{11}+Ric(T)_{22}+2Ric(T)_{33}=(Ric(T)_{11}+Ric(T)_{33})+(Ric(T)_{22}+Ric(T)_{33})\geq-\frac{4\gamma_{max}}{n}\scal(T),
\end{equation*}
we obtain

\begin{equation*}
    \begin{split}
        RHS\geq& -\frac{1+(n-2)b_{max}}{1+(n-2)b}\sqrt{2a}\frac{n^2-7n+14}{n^2-5n+4}(Ric(T)_{11}+Ric(T)_{22}-2T_{1212})\\
        &+\frac{b_{max}-b}{1+(n-2)b}((1+\lambda^2)Ric(T)_{11}+(1+\lambda^2)Ric(T)_{22}+2Ric(T)_{33}+2\lambda^2Ric(T)_{44})\\
        &+\frac{4}{n}\left(\frac{a_{max}-a}{1+2(n-1)a}-\frac{b_{max}-b}{1+(n-22)b}\right)\scal(T)\\
        &+\frac{1+(n-2)b_{max}}{1+(n-2)b}\sqrt{2a}(Ric(T)_{11}+Ric(T)_{22})\\
        &+\frac{2}{n}\left(\frac{2(n-1)(a_{max}-a)}{1+2(n-11)a}-\frac{(n-2)(b_{max}-b)}{1+(n-2)b}\right)\sqrt{2a}\scal(T).
    \end{split}
\end{equation*}

Notice that we can give a lower bound for the term $T_{1212}$ in terms of $Ric(T)_{11}+Ric(T)_{22}$ and $\scal(T)$, as follows:
\begin{equation*}
    \begin{split}
        0\leq &\sum_{3\leq p \neq q \leq n}(T_{1212}+T_{1p1p}+T_{2q2q}+T_{pqpq})\\
        =& (n^2-7n+14)T_{1212}+(n-5)(Ric(T)_{11}+Ric(T)_{22})+\scal(T).
    \end{split}
\end{equation*}

This gives 
\begin{equation*}
    \begin{split}
        RHS\geq&-\frac{1+(n-2)b_{max}}{1+(n-2)b}\frac{2\sqrt{2a}}{n^2-5n+4}\scal(T)\\
        +&\frac{4}{n}\left(\frac{a_{max}-a}{1+2(n-1)a}-(1+\gamma_{max})\frac{b_{max}-b}{1+(n-2)b}\right)\scal(T)\\
        +&\frac{2}{n}\left(\frac{2(n-1)(a_{max}-a)}{1+2(n-1)a}-\frac{(n-2)(b_{max})-b}{1+(n-2)b}\right)\sqrt{2a}\scal(T).
    \end{split}
\end{equation*}

Thus, $RHS\geq0$ by Lemma \ref{lemma connect cone inequalities}. $\hfill\square$

\hspace*{\fill}

Now we are left to show the following theorem:

\begin{theorem}\label{thm charactrization of second cone}(cf. Theorem 4.6 of \cite{Brendle19})
    For each $0<b\leq b_{max}$, the cone $\Tilde{C}(b)$ is transversally invariant under the Hamiltonian ODE $\frac{\mathrm{d}}{\mathrm{d}t}R=Q(R)$.
\end{theorem}

Here, we need the following two computation results in \cite{Brendle19}.

\begin{lemma}\label{lemma BrendlelemmaA1}(Lemma A. 1 of \cite{Brendle19}) Let $0\leq \zeta\leq 1$ and $0<\rho\leq 1$. Assume that $H$ is a symmetric bilinear form with the property that the largest eigenvalue of $H$ is bounded from above by $\frac{1}{2}\mathrm{tr}(H)$ and the sum of the two smallest eigenvalues of $H$ is bounded from below by $\frac{2(1-\zeta)}{n}\mathrm{tr}(H)$. Then

\begin{equation*}
    \begin{split}
        \frac{n-2}{n}\mathrm{tr}(H)(H_{11}+H_{22})&-\rho ((H_0^2)_{11}+(H_0^2)_{22})\\
        \geq&\frac{2}{n^2}\left((n-2)(1-\zeta)-2\zeta^2\rho\frac{n^2-2n+2}{(n-2)^2}\right)\mathrm{tr}(H)^2
    \end{split}
\end{equation*}
for every pair of orhonormal vectors $\{e_1,e_2\}$, where $H_0$ denotes the tracefree part of $H$.

\end{lemma}

\begin{lemma}\label{lemmaBrendleA8}(Lemma A. 8 of \cite{Brendle19})
    Let $S$ be an algebraic curvaure tensor, and let $H$ be a symmetric bilinear form on $\mathbb{R}^n$ such that
    \begin{equation*}    Z:=S_{1313}+\lambda^2S_{1414}+S_{2323}+\lambda^2S_{2424}-2\lambda S_{1234}+(1-\lambda^2)(H_{11}+H_{22})\geq0
    \end{equation*}
for every orthonormal four-frame $\{e_1,e_2,e_3,e_4\}$ and every $\lambda\in [0,1]$. Moreover, suppose that $Z=0$ for one particular orthonormal four-frame $\{e_1,e_2,e_3,e_4\}$ and one particular $\lambda\in [0,1)$. Then, for this particular orthonormal four-frame $\{e_1,e_2,e_3,e_4\}$ and this $\lambda\in [0,1)$, we have

\begin{equation*}
    \begin{split}
Q(S)_{1313}&+\lambda^2Q(S)_{1414}+Q(S)_{2323}+\lambda^2Q(S)_{2424}-2\lambda Q(S)_{1234}\\
        &+(H\wedge H)_{1313}+\lambda^2(H\wedge H)_{1414}+(H\wedge H)_{2323}+\lambda^2(H\wedge H)_{2424}-2\lambda(H\wedge H)_{1234}\\
        &+2(1-\lambda^2)((S\star H)_{11}+(S\star H)_{22})\\
        \geq&(1+\lambda^2)(H_{11}+H_{22})^2,
    \end{split}
\end{equation*}
where we define $(S\star H)_{ik}=\sum_{p,q=1}^nS_{ipkq}H_{pq}$.
\end{lemma}

\noindent\textit{Proof of Theorem \ref{thm charactrization of second cone}.} It suffices to show that $\Tilde{\mathcal{E}}(b)$ is transversally invariant under $\frac{\mathrm{d}}{\mathrm{d}t}S=Q(S)+D_{a,b}(S)$ for each $0<b\leq\Tilde{b}_{max}$. Suppose $S\in \Tilde{\mathcal{E}}(b)\backslash\{0\}$, then
\begin{equation*}
    Z:=S_{1313}+\lambda^2S_{1414}+S_{2323}+\lambda^2S_{2424}-2\lambda S_{1234}+\sqrt{2a}(1-\lambda^2)(Ric(S)_{11}+Ric(S)_{22})\geq0
\end{equation*}
for each orthonormal four-frame $\{e_1,e_2,e_3,e_4\}$ and every $\lambda\in [0,1]$.

To show the desired result, we assume $Z=0$ for some orthonormal four-frame $\{e_1,e_2,e_3,e_4\}$ and some $\lambda\in [0,1]$. We distinguish two cases:

\noindent\textit{Case 1.} Suppose that $\lambda\in [0,1)$. Using the formula for $D_{a,b}(S)$ in Section 2, we deduce that

\begin{equation*}
    \begin{split}
        \frac{\mathrm{d}}{\mathrm{d}t}Z\geq& Q(S)_{1313}+\lambda^2Q(S)_{1414}+Q(S)_{2323}+\lambda^2Q(S)_{2424}-2\lambda Q(S)_{1234}\\
        &+2a(Ric(S)\wedge Ric(S))_{1313}+2a\lambda^2(Ric(S)\wedge Ric(S))_{1414}\\
        &+2a(Ric(S)\wedge Ric(S))_{2323}+2a\lambda^2(Ric(S)\wedge Ric(S))_{2424}-4a\lambda(Ric(S)\wedge Ric(S))_{1234}\\
        &+2\sqrt{2a}(1-\lambda^2)((S\star Ric(S))_{11}+(S\star Ric(S))_{22})\\
        &+\frac{4(n-2)a}{n}\sqrt{2a}(1-\lambda^2)\scal(S)(Ric(S)_{11}+Ric(S)_{22})\\
        &-4b\sqrt{2a}(1-\lambda^2)((Ric_0(S)^2)_{11}+(Ric_0(S)^2)_{22})+\frac{8a}{n^2}\sqrt{2a}(1-\lambda^2)\scal(S)^2.
    \end{split}
\end{equation*}

Applying Lemma \ref{lemmaBrendleA8} with $H=\sqrt{2a}Ric(S)$, we obtain 

\begin{equation*}
    \begin{split}
        Q(S)_{1313}&+\lambda^2Q(S)_{1414}+Q(S)_{2323}+\lambda^2Q(S)_{2424}-2\lambda Q(S)_{1234}\\
        &+2a(Ric(S)\wedge Ric(S))_{1313}+2a\lambda^2(Ric(S)\wedge Ric(S))_{1414}\\
        &+2a(Ric(S)\wedge Ric(S))_{2323}+2a\lambda^2(Ric(S)\wedge Ric(S))_{2424}-4a\lambda(Ric(S)\wedge Ric(S))_{1234}\\
        &+2\sqrt{2a}(1-\lambda^2)((S\star Ric(S))_{11}+(S\star Ric(S))_{22})\geq0.
    \end{split}
\end{equation*}

Since $S\in \Tilde{\mathcal{E}}(b)\subset C_{PIC}$, we have that the largest eigenvalue of $Ric(S)$ is bounded from above by $\frac{1}{2}\scal(S)$. Moreover, by Lemma \ref{lemma elemenary2 second cone invariant}  , the sum of the two smallest eigenvalues of $Ric(S)$ is bounded from below by $\frac{2(1-\zeta)}{n}\scal(S)$. Applying Lemma \ref{lemma BrendlelemmaA1} for $H=Ric(S)$ and $\rho =\frac{b}{a}$ , we obtain

\begin{equation*}
    \begin{split}
        \frac{(n-2)a}{n}\scal(S)(Ric(S)_{11}+Ric(S)_{22})&-b((Ric_0(S)^2)_{11}+(Ric_0(S)^2)_{22})\\
        \geq&\frac{2}{n^2}\left(a(n-2)(1-\zeta)-2b\zeta^2\frac{n^2-2n+2}{(n-2)^2}\right)(\scal(S))^2.
    \end{split}
\end{equation*}

Putting these facts together, and notice that $(1+(n-2)(1-\zeta))>2\zeta^2\frac{n^2-2n+2}{(n-2)^2}$ by Lemma \ref{lemma elementary second cone invariant} , we conclude that 
\begin{equation*}
    \frac{\mathrm{d}}{\mathrm{d}t}Z\geq\frac{8\sqrt{2a}(1-\lambda^2)}{n^2}\left(a(1+(n-2)(1-\zeta))-2b\zeta^2\frac{n^2-2n+2}{(n-2)^2}\right)(\scal(S))^2>0
\end{equation*}

as $\lambda\in [0,1)$.

\noindent\textit{Case 2.} Suppose now that $\lambda=1$,then we calculate that

\begin{equation*}
    \begin{split}
        \frac{\mathrm{d}}{\mathrm{d}t}Z=&Q(S)_{1313}+Q(S)_{1414}+Q(S)_{2323}+Q(S)_{2424}-2Q(S)_{1234}\\
        &+2a[(Ric(S)\wedge Ric(S))_{1313}+(Ric(S)\wedge Ric(S))_{1414}\\
        &+(Ric(S)\wedge Ric(S))_{2323}+(Ric(S)\wedge Ric(S))_{2424}-2(Ric(S)\wedge Ric(S))_{1234}]\\
    \end{split}
\end{equation*}

Since $S\in C_{PIC}$ and $Z=S_{1313}+S_{1414}+S_{2323}+S_{2424}-2Q(S)_{1234}\geq0$, noticing that the PIC condition is invariant under Hamiltonian ODE $\frac{\mathrm{d}}{\mathrm{d}t}S=Q(S)$, we obtain 

\begin{equation*}
    Q(S)_{1313}+Q(S)_{1414}+Q(S)_{2323}+Q(S)_{2424}-2Q(S)_{1234}\geq0.
\end{equation*}

Now it suffices to show that the algebraic curvature operator $Ric(S)\wedge Ric(S)$ lies in the interior of the cone $ C_{PIC}$. To this end, let $\zeta,\eta\in \mathbb{C}^n$ be linearly independent complex vectors satisfying $g(\zeta,\zeta)=g(\zeta,\eta)=g(\eta,\eta)=0$, it suffices to show that $(Ric(S)\wedge Ric(S))(\zeta,\eta,\Bar{\zeta},\Bar{\eta})>0$.  As $Ric(S)$ is unitarily diagonalizable on $\mathrm{span}_{\mathbb{C}}\{\eta,\zeta\}$, we can find $z,w\in\mathrm{span}_{\mathbb{C}}\{\eta,\zeta\}$ such that $g(z,\Bar{z})=g(w,\Bar{w})=2, g(z,\Bar{w})=0$ and $Ric(S)(z,\Bar{w})=0$. The identities $g(\zeta,\zeta)=g(\zeta,\eta)=g(\eta,\eta)=0$ give $g(z,z)=g(w,w)=g(z,w)=0$. Consequently, we may write $z=e_1+ie_2$ and $w=e_3+ie_4$ for some orthonormal four-frame $\{e_1,e_2,e_3,e_4\}\subset \mathbb{R}^n$. Using the identity $Ric(S)(z,\Bar{w})=0$, and by Lemma  we have that the sum of two smallest eigenvalues of $Ric(S)$ is strictly positive, we obtain

\begin{equation*}
    (Ric(S)\wedge Ric(S))(z,w,\Bar{z},\Bar{w})=2(Ric(S)_{11}+Ric(S)_{22})(Ric(S)_{33}+Ric(S)_{44})>0.
\end{equation*}

Since $\mathrm{span}_{\mathbb{C}}\{\eta,\zeta\}=\mathrm{span}_{\mathbb{C}}\{z,w\}$, we have $(Ric(S)\wedge Ric(S))(\zeta,\eta,\Bar{\zeta},\Bar{\eta})>0$. Combining these facts together, we deduce that $\frac{\mathrm{d}}{\mathrm{d}t}Z>0$. $\hfill\square$

\hspace*{\fill}

Arriving here, we are finally able to prove Theorem \ref{thm main cone}.

\noindent\textit{Proof of Theorem \ref{thm main cone}.} Since $\Tilde{C}(\Tilde{b}_{max})=C(b_{max})$, let $B=b_{max}+\Tilde{b}_{max}$, define
\begin{equation*}
    \hat{C}(b)=\begin{cases}
        C(b), & b\in (0,b_{max}]\\
        \Tilde{C}(B-b), & b\in(b_{max},B).
    \end{cases}
\end{equation*}
Then $\hat{C}(b)$ is a family of closed, convex, $O(n)$-invariant cones that depends continuously on $b\in (0,B)$. By Theorem \ref{thm goal for first family} and Theorem \ref{thm charactrization of second cone}, for each $b\in (0,B)$, the cone $\hat{C}(b)$ is transversally invariant under the Hamilton ODE.

Moreover, by Definition \ref{def firstcone}, the properties of the family $C(b)$ guarantee that there exists a small $\beta_0>0$ such that  \begin{equation*}
            \begin{split}
                (\{R:R-\theta scal\mathrm{id}\wedge\mathrm{id}\in C_{PIC}\}&\cap\{R: Ric_{11}+Ric_{22}-\theta scal+N\geq0\} )\\
                &\subset\{R:R+N\mathrm{id}\wedge
                \mathrm{id}\in \hat{C}(\beta_0)\}.
            \end{split}
        \end{equation*}

Now, let $\{\beta_j\}$ be an increasing sequence such that $\lim_{j\rightarrow +\infty}\beta_j=B$. By Definition \ref{def second family}, there exists a sequence of positive reals $\{\varepsilon_j\}$ with $\lim_{j\rightarrow +\infty}\varepsilon_j=0$, and
        \begin{equation*}
            \hat{C}(\beta_j)\subset \{R:R+\varepsilon_jscal\mathrm{id}\wedge\mathrm{id}\in C_{PIC1}\}.
        \end{equation*}

This proves Theorem \ref{thm main cone}.$\hfill\square$

\bibliographystyle{alpha}
\bibliography{PICbib}

\end{document}